\documentclass[10pt]{article}
\usepackage{amsmath,amssymb,amsthm, tikz,tikz-cd, subcaption, float}
\usepackage{clock}
\usepackage[clock]{ifsym}
\usepackage{listings}
\usetikzlibrary{arrows.meta}
\usetikzlibrary{positioning}
\usetikzlibrary{calc}
\usetikzlibrary{decorations.markings}
\usetikzlibrary{arrows}

\usepackage{amsmath,amssymb,amsfonts,amsthm,mathtools,mathrsfs}
\usepackage{cancel,bm,array}
\usepackage[margin=1in]{geometry}
\usepackage{tikz}
\usetikzlibrary{arrows.meta,positioning,calc,decorations.markings,arrows}
\tikzset{> /.tip = {Latex[round,length=.3cm]}}
\tikzset{every path/.style={thin}}
\tikzset{baseline={([yshift=-.5ex]current bounding box.center)}}
\tikzstyle{mid>}=[decoration={markings, mark=at position 0.5 with {\arrow{>}}}, postaction={decorate}]
\tikzstyle{mid<}=[decoration={markings, mark=at position 0.5 with {\arrowreversed{>}}}, postaction={decorate}]
\tikzstyle over=[preaction={draw,line width=6pt,white}]
\usepackage{xcolor}
\definecolor{ForestGreen}{RGB}{34,139,34}

\newtheoremstyle{spacedplain}
  {12pt}
  {12pt}
  {\itshape}
  {}
  {\bfseries}
  {.}
  {.5em}
  {}
\newtheoremstyle{spaceddefn}%
  {12pt}{12pt}
   {\itshape}
  {}{\bfseries}{.}{.5em}{}
\theoremstyle{spacedplain}
\newtheorem{theorem}{Theorem}
\numberwithin{theorem}{section}
\newtheorem{proposition}[theorem]{Proposition}
\newtheorem{lemma}[theorem]{Lemma}
\newtheorem{corollary}[theorem]{Corollary}
\theoremstyle{spaceddefn}
\newtheorem{definition}[theorem]{Definition}
\newtheorem{example}[theorem]{Example}
\newtheorem{remark}[theorem]{Remark}

\newcommand{\excise}[1]{}                 

\newcommand{\Ww}{\mathbb{W}}              

\newcommand{\Zz}{\mathbb{Z}}
\newcommand{\Gg}{\mathbb{G}}

\newcommand{\inv}{\mathfrak{sl}_3}
\def\la{\langle}
\def\ra{\rangle}
\def\llangle{\la\!\la}
\def\rrangle{\ra\!\ra}

\newcommand{\leftlangle}{\left\la\!\left\la}
\newcommand{\rightrangle}{\right\ra\!\right\ra}
\usepackage[colorlinks=true,linkcolor=blue,citecolor=blue,urlcolor=blue,%
            breaklinks=true,bookmarks=true]{hyperref}

\newcommand{\Kp}{K_{p+1}}
\newcommand{\lam}{\lambda}

\usepackage{hyperref}

\tikzset{> /.tip = {Latex[round,length=.3cm]}}
\tikzset{every path/.style={very thick}}
\tikzset{baseline={([yshift=-.5ex]current bounding box.center)}}
\tikzstyle{mid>}=[decoration={markings, mark=at position 0.5 with {\arrow{>}}}, postaction={decorate}]
\tikzstyle{mid<}=[decoration={markings, mark=at position 0.5 with {\arrowreversed{>}}}, postaction={decorate}]
\tikzstyle over=[preaction={draw,line width=6pt,white}]

\newcommand{\idL}{%
  \begin{tikzpicture}[scale = 1]
			\draw (0,0) to (0,1);
            \draw (0.25,0) to (0.25,1);
            \draw (0.5,0) to (0.5,1);
\end{tikzpicture}%
}
\newcommand{\rw}{%
  \begin{tikzpicture}[scale = 1]
            \draw (-0.25,0) to (-0.25, 1);
            \draw (0,0) to (.25,.25);
            \draw (0.5,0) to (.25,.25);	
            \draw (.25,.25) to (.25,.75);
            \draw (.25,.75) to (0,1);
            \draw (.25,.75) to (0.5,1);		
\end{tikzpicture}%
}
\newcommand{\lw}{%
  \begin{tikzpicture}[scale = 1]
            \draw (0.75,0) to (0.75, 1);
            \draw (0,0) to (.25,.25);
            \draw (0.5,0) to (.25,.25);	
            \draw (.25,.25) to (.25,.75);
            \draw (.25,.75) to (0,1);
            \draw (.25,.75) to (0.5,1);	
\end{tikzpicture}%
}
\newcommand{\rl}{%
  \begin{tikzpicture}[scale = 0.5]
    \draw (-0.5,0) to (-0.5, 1);
    \draw (0,0) to (.25,.25);
    \draw (0.5,0) to (.25,.25);	
    \draw (.25,.25) to (.25,.75);
    \draw (.25,.75) to (0,1);
    \draw (.25,.75) to (0.5,1);	
\begin{scope}[shift={(-0.50,1)}]
    \draw (1,0) to (1, 1);
    \draw (0,0) to (.25,.25);
    \draw (0.5,0) to (.25,.25);	
    \draw (.25,.25) to (.25,.75);
    \draw (.25,.75) to (0,1);
    \draw (.25,.75) to (0.5,1);
\end{scope}            
\end{tikzpicture}%
}
\newcommand{\lr}{%
  \begin{tikzpicture}[scale = 0.5]
    \draw (1,0) to (1, 1);
    \draw (0,0) to (.25,.25);
    \draw (0.5,0) to (.25,.25);	
    \draw (.25,.25) to (.25,.75);
    \draw (.25,.75) to (0,1);
    \draw (.25,.75) to (0.5,1);
\begin{scope}[shift={(0.50,1)}]
    \draw (-0.5,0) to (-0.5, 1);
    \draw (0,0) to (.25,.25);
    \draw (0.5,0) to (.25,.25);	
    \draw (.25,.25) to (.25,.75);
    \draw (.25,.75) to (0,1);
    \draw (.25,.75) to (0.5,1);
\end{scope}
\end{tikzpicture}%
}
\newcommand{\cc}{%
  \begin{tikzpicture}[scale=0.5]
    \draw (-0.5,0.5) to[out=-90, in=-90] (0.5,0.5);
    \draw (0.5,0.5) to (0.5,1);
    \draw (-0.5,0.5) to (-0.5,1);
    \draw (0,0.2) to (0, 1);
    
    \draw (-0.5,-0.5) to[out=90, in=90] (0.5,-0.5);
    \draw (0.5,-0.5) to (0.5,-1);
    \draw (-0.5,-0.5) to (-0.5,-1);
    \draw (0,-0.2) to (0, -1);
\end{tikzpicture}%
}

\newcommand{\state}[3]{
  \draw[mid>] ( 0.25,1.5) to[out=0,in=0]     ( 0.25,0);
  \draw[mid<] (-0.25,0)   to[out=180,in=180] (-0.25,1.5);
  \csname #1tile\endcsname{0}
  \csname #2tile\endcsname{0.5}
  \csname #3tile\endcsname{1}}

\newlength{\cellpad}
\newcommand{\vc}[1]{\ensuremath{\vcenter{\kern\cellpad\hbox{#1}\kern\cellpad}}}
 
\newcommand{\ID}{\vc{\idL}}
\newcommand{\RW}{\vc{\rw}}   
\newcommand{\LW}{\vc{\lw}}   
\newcommand{\LR}{\vc{\lr}}   
\newcommand{\RL}{\vc{\rl}}   
\newcommand{\CC}{\vc{\cc}}   
 
\newcommand{\webalg}{\mathrm{Web}_3}

\title{Combinatorial characterizations of the quantum $\mathfrak{sl}_3$ invariant of positive alternating links and 3-braids}
\author{Matthew Harper, Joshua O'Connor, Jos\'e S\'anchez Mench\'en,  Angie Wu}
\date{November 22, 2026}

\begin{document}
\maketitle

\begin{abstract}
\noindent The $\mathfrak{sl}_3$ link polynomial is a higher rank relative of the Jones polynomial which can be studied entirely diagrammatically using Kuperberg's web calculus. Building on the work of Harper and Kalfagianni, we construct formulae for the fourth and fifth coefficients of the quantum $\mathfrak{sl}_3$ link invariant on positive alternating link diagrams in terms of data read from the associated Seifert graph. Then, we give a characterization of the invariant on positive alternating fibered links in terms of elementary symmetric polynomials. We also discuss behavior of our formulae under connected sums. We give formulations for the $\mathfrak{sl}_3$ invariant of even powers of the Garside element $\Delta_3\in B_3$ and a closed formula for the invariant on positive $3$-braids, resulting in a formula for the invariant on all $3$-braids.
\end{abstract}

\section{Introduction}
The Jones polynomial is the simplest and most extensively studied member of a family of quantum link invariants obtained by applying the Reshetikhin--Turaev construction to the quantized enveloping algebras $U_q(\mathfrak{g})$. The case $\mathfrak{g} = \mathfrak{sl}_2$ recovers the Jones polynomial, and the next case $\mathfrak{g} = \mathfrak{sl}_3$ produces the {quantum $\inv$ invariant} $\llangle L \rrangle \in \Zz[q,q^{-1}]$. Kuperberg \cite{Kuperberg} showed that this invariant admits a purely combinatorial description through a web calculus.\\

In order to relate these link polynomials to the topology and geometry of knots, the first named author and Kalfagianni~\cite{HK} studied the first three leading coefficients of $\llangle L \rrangle$ for positive links. They expressed these coefficients in terms of quantities which could be read directly from the link diagram, or from its Seifert graph. They showed, among other things, that the invariant detects fiberedness and obstructs the presentation of a positive link as a positive braid. Their work also uses ideas of \cite{DasbachLin} which considers the Jones polynomial.\\

In this paper, we build directly on the work of Harper and Kalfagianni, by studying the fourth and fifth coefficients of the $\inv$ polynomial for positive alternating links. Recall that an oriented link diagram is \emph{alternating} if when traversing each component of the diagram along its orientation, one alternates between over-crossings and under-crossings. An oriented link diagram is \emph{positive} if all of its crossings are positive. A link is alternating (resp.\! positive) if it admits an alternating (resp.\! positive) diagram. By work of \cite{nak}, a link which admits a positive diagram and an alternating diagram also admits a diagram which is both positive and alternating, see also \cite{Stoimenow}. We maintain the positivity assumption of \cite{HK} in order to sort state sum contributions of the polynomial by degree. The assumption that we consider only alternating diagrams is to avoid introducing {mixed states}, seen in \cite[{Section 5.1}]{HK} which occur at type-2 Seifert circles and would increase the complexity of our computations for the fourth and fifth coefficients.\\

To describe our results, let $L$ be a positive link and write
\[
\llangle L \rrangle = \sum_{i\geq 1} \gamma_i q^{n-2(i-1)}.
\]
It was shown in \cite{HK} that $\llangle L \rrangle$ indeed takes this form and $n=v-e$, where $v$ and $e$ are the number of vertices and edges in the Seifert graph $\Gg_W$. We also write $\Gg_W'$ for the reduced graph of $\Gg_W$, obtained by removing all multiple edges between all pairs of vertices, and denote the number of reduced edges by $e'$.

\begin{theorem}
For any positive alternating link:
\begin{align*}
    \gamma_4 &= \binom{v-e'+2}{3} + (\mu_1-1)\,(v-e') + e' - \mu_2\,,
    \\
    \gamma_5 &= \binom{v-e'+3}{4} + (\mu_1-2)\binom{v-e'+1}{2} + (e'-\mu_2+1)(v-e') + \binom{\mu_1}{2} - e' + \mu_3 + \sigma.
\end{align*}
Here $\mu_i$ is the number of reduced edges with multiplicity different from $i$ and $\sigma$ is the number of four-cycles in $\Gg_W'$.
\end{theorem}
In the notation of \cite{HK}, $\mu$ is equivalent to our definition of $\mu_1$.\\

These calculation motivated the investigation of additional terms. If we assume that the Seifert graph is a tree, which corresponds to the fibered case, we give a closed formula for the $\mathfrak{sl}_3$ invariant in terms of symmetric polynomials. We use this to show the invariant can be realized as a Tutte evaluation.

\begin{theorem}
    For every link $L$ with a positive alternating connected link diagram, whose reduced Seifert graph $\Gg_W'$ is a tree, we have
\begin{equation*}
    \llangle L \rrangle = q^{-2e}\sum_{r=0}^{e'}[3]^{\,v-r}\,e_r(u_1,\dots,u_{e'}),
  \qquad u_i = (-q^{-2})^{z_i}-1,
\end{equation*}
where $e_r$ is the $r$-th elementary symmetric polynomial and $z_i$ is the multiplicity of the $i$-th reduced edge.
\end{theorem}

\begin{theorem}
    Let $L$ be a connected link which admits a positive alternating diagram, and whose reduced Seifert graph $\Gg_W'$ is a tree. Then \[\llangle L\rrangle = q^{-2e}[3](-q^{-2}-1)^{e'}T(M(\Gg_W); \frac{-q^4}{q^{2}+1}, -q^{-2}).\]
\end{theorem}

We consider a second combinatorial characterization of the $\inv$ polynomial, now on arbitrary braids on three strands, in $B_3$. Our approach uses the Garside element $\Delta_3$ and the algebra of three strand webs $\webalg$. The algebra is equipped with a basis of non-elliptic vectors 
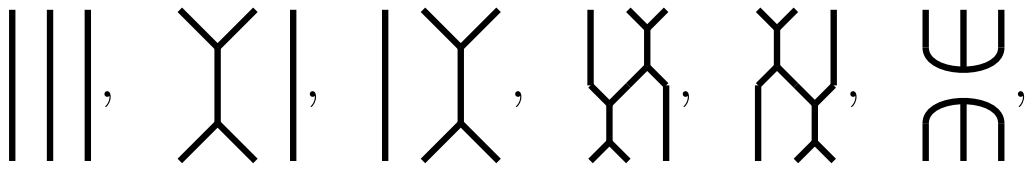
\[\quad\idL\,, \quad\lw\,,\quad\rw\, ,\quad\lr\,,\quad\rl\,,  \quad\cc\,,\]
and we describe its complete multiplication table in Figure \ref{fig:mult}. \\

Recall that $\Delta_3^2$ is central and that for any braid $b\in B_3$ there exists a positive braid $p$ and integer $m=2k$ such that $b=p\cdot \Delta_3^{2k}$. After giving a closed form expression for the components of any power of $\Delta_3$ in the basis, we use this result to compute the invariant of the braid closure $\llangle \hat{b}\rrangle$ in terms of the data of its twist regions. A similar approach was used for computing the Khovanov homology of 3-braids in \cite{VGSGarside}.

\begin{theorem}
Take $b\in B_3$, with $b=p\cdot \Delta_3^{-2k}$ for some $k\in \mathbb Z_{\ge 0}$ and $p=\sigma_1^{n_1}\cdots \sigma_2^{n_d}$ up to Markov moves. Then, 
\begin{equation*}
    \begin{aligned}
        \llangle \hat{b}\rrangle &= q^{12k-2\mathbf{n}}\sum_{k_1, \cdots, k_d}\binom{n_1}{k_1}\cdots \binom{n_d}{k_d}((-q)^{-\mathbf{k}})[2]^{\sum_{k_i:k_i\ge 2}(k_i-1)}
        \\&\qquad
        \cdot 
        \bigg(\llangle \mathbb{W}(T'(0))\rrangle +q^{3k}\{2\}[2][k]_{q^3}\mathbb{W}(T'(1))\rrangle+q^{3k}\{2\}[2][k]_{q^3}\llangle \mathbb{W}(T'(2))\rrangle
        \\&\qquad\qquad -q^{3k}\{2\}[k]_{q^3}\mathbb{W}(T'(3))\rrangle
        -q^{3k}\{2\}[k]_{q^3}\mathbb{W}(T'(4))\rrangle+q^{6k}\{2\}\{3\}[k]_{q^3} [2]^{\mathbf{t}+1}[3]\bigg) 
    \end{aligned}
\end{equation*}
where $T'(i)$ is a cyclic binary sequence associated to the twist regions in a given state augmented by a non-elliptic basis vector.
\end{theorem}

The outline of our paper is as follows.
\begin{itemize}
    \item \textbf{Section 2}: An introduction to the $\mathfrak{sl}_3$ invariant and webs.
    \item \textbf{Section 3}: Computation of the fourth and fifth coefficients, and a formula for the invariant when the reduced Seifert graph is a tree. 
    \item \textbf{Section 4}: Applications of formulae to connected sums.
    \item \textbf{Section 5}: A characterization of the $\mathfrak{sl}_3$ invariant on $3$-braids.
\end{itemize}

\subsection*{Tool Use} AI was used to help generate the code at \cite{AngieCode} by describing the steps of the formula. It was also used to convert the references of the original draft to bibtex format. 

\subsection*{Acknowledgements}
We would like to thank Bob Bell for running the 2026 SURIEM REU, where most of this work took place, and for making this collaboration possible. The authors thank Efstratia Kalfagianni and Anthony Stone for helpful conversations. MH was partially supported through the NSF-RTG grant \#DMS-2135960. The remaining authors were supported by NSF Grant \#2244461 through the SUREIM REU at Michigan State University.  MH and JO thank the Michigan State University SROP program for additional support. 

\section{The quantum \texorpdfstring{$\mathfrak{sl}_3$}{sl3} invariant}\label{sec:defofinv}

The quantum $\inv$ invariant is defined diagrammatically. Given an oriented link diagram $D$, we resolve each crossing according to the following relations:
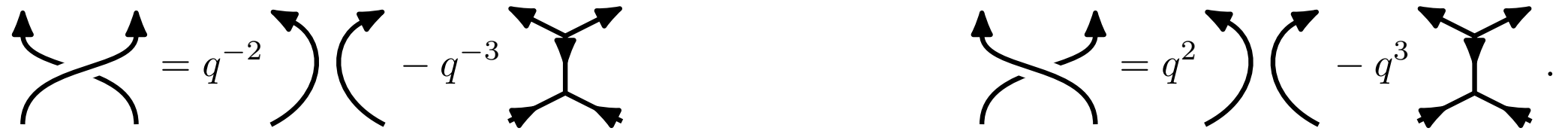
\begin{align}\label{eq.crossingrelations}
        \begin{tikzpicture}[scale = 1]
            \draw[->] (1,0) to [out=90,in=-90] (0,1);
            \draw[over,->] (0,0) to[out=90,in=-90] (1,1);
        \end{tikzpicture}
        =q^{-2}
        \begin{tikzpicture}[scale = 1]
            \draw[->] (0,0) .. controls (.5,.25) and (.5,.75) .. (0,1);
            \draw[->] (1,0) .. controls (.5,.25) and (.5,.75) .. (1,1);
        \end{tikzpicture}-
        q^{-3}
        \begin{tikzpicture}[scale = 1] 
            \draw[mid>] (0,0) to (.5,.25);
            \draw[mid>] (1,0) to (.5,.25);	
            \draw[mid<] (.5,.25) to (.5,.75);
            \draw[->] (.5,.75) to (0,1);
            \draw[->] (.5,.75) to (1,1);		
        \end{tikzpicture}\hspace{3cm}
        \begin{tikzpicture}[scale = 1]
            \draw[->] (0,0) to [out=90,in=-90] (1,1);
            \draw[over,->] (1,0) to[out=90,in=-90] (0,1);
        \end{tikzpicture}
        =q^{2}
        \begin{tikzpicture}[scale = 1]
            \draw[->] (0,0) .. controls (.5,.25) and (.5,.75) .. (0,1);
            \draw[->] (1,0) .. controls (.5,.25) and (.5,.75) .. (1,1);
        \end{tikzpicture}-
        q^{3}
        \begin{tikzpicture}[scale = 1]
            \draw[mid>] (0,0) to (.5,.25);
            \draw[mid>] (1,0) to (.5,.25);	
            \draw[mid<] (.5,.25) to (.5,.75);
            \draw[->] (.5,.75) to (0,1);
            \draw[->] (.5,.75) to (1,1);		
        \end{tikzpicture}\,.
\end{align}

By replacing each crossing in this manner, we introduce a $\Zz[q,q^{-1}]$-linear combination of web diagrams, trivalent graphs embedded in $S^2$ where each vertex is either a source or a sink.. Once every crossing is resolved we may use the following \emph{square}, \emph{bubble}, and \emph{circle} relations to simplify the resulting diagrams further:
   
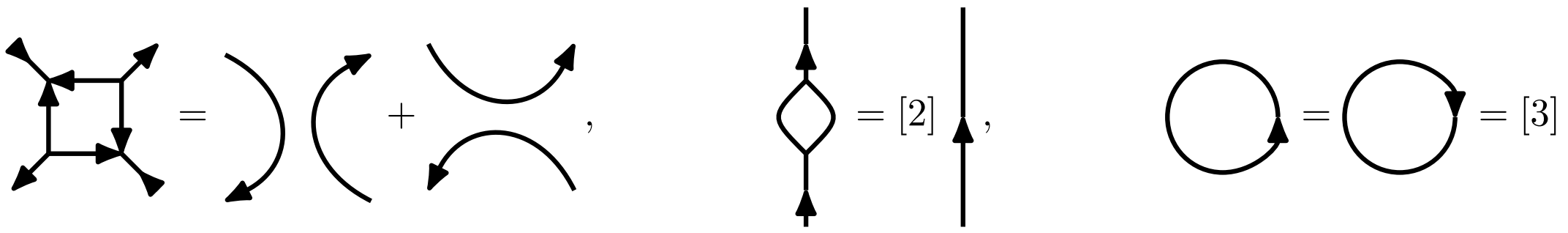
\begin{align}\label{eq.webrelations}
        \begin{tikzpicture}[scale=0.66]
			\draw[->] (0,0) to (0,1);
			\draw[->] (1,1) to (1,0);
			\draw[->] (0,0) to (1,0);
			\draw[->] (1,1) to (0,1);
			\draw[->] (0,0) to (-.5,-.5);
			\draw[mid>] (1.5,-.5) to (1,0);
			\draw[mid>] (-.5,1.5) to (0,1);
			\draw[->] (1,1) to (1.5,1.5);
		\end{tikzpicture}
		=
		\begin{tikzpicture}[scale = 0.66,baseline=10]
			\draw[->] (-.5,1.5) .. controls (.5, 1) and (.5,0) ..  (-.5,-.5);
			\draw[,<-] (1.5,1.5) .. controls (.5, 1) and (.5,0) ..  (1.5,-.5);
		\end{tikzpicture}
		+
		\begin{tikzpicture}[scale = 0.66]
			\draw[->] (-.5,1.5) .. controls (0, .5) and (1,.5) ..  (1.5,1.5);
			\draw[,<-] (-.5,-.5) .. controls (0,.5) and (1,.5) ..  (1.5,-.5);
		\end{tikzpicture}
        \,,\hspace{1.5cm}
        \begin{tikzpicture}[scale = 0.66]
			\draw[->] (0,0) to (0,.5);
			\draw[] (0,.5) to (0,1);
			\draw[] (0,1) .. controls (-.5,1.5) .. (0,2);
			\draw[] (0,1) .. controls (.5,1.5) .. (0,2);
			\draw[->] (0,2) to (0,2.5);
			\draw[] (0,2.5) to (0,3);		
		\end{tikzpicture} = [2]\,\,
        \begin{tikzpicture}[scale = 0.66]
			\draw[->] (0,0) to (0,1.5);
			\draw[] (0,1.5) to (0,3);
		\end{tikzpicture}\,
        ,\hspace{1.5cm}
		\begin{tikzpicture}[scale = 0.5]
			\draw[->] (0,0) arc(0:360:1);
		\end{tikzpicture} =
		\begin{tikzpicture}[scale = 0.5]
			\draw[<-] (0,0) arc(0:360:1);
		\end{tikzpicture}
        = [3]
   \end{align}
where $[2] = q+q^{-1}$ and $[3] = q^2 + 1 +q^{-2}$. More generally, we write for $k\in\mathbb{Z}_{\geq0}$,
\begin{align}
    \label{eq.qnumber}
    [k]_q=\frac{q^k-q^{-k}}{q-q}=q^{k-1}+q^{k-3}+\cdots+q^{-(k-3)}+q^{-(k-1)}\qquad \mbox{with $[0]_q=0$.}
\end{align}
We will omit the subscript and write $[k]$ to mean $[k]_q$.

\begin{lemma}[{\cite[Lemma B.1]{Ohtsuki}}] \label{bigonsquare}
	Let $\Ww$ be a nonempty web embedded on a 2-sphere $S^2$. 	 Then, at least one of the regions of $S^2\setminus \Ww$ is a bigon or a square. Therefore, any closed web can be evaluated to a multiple of the empty diagram.
\end{lemma}

Given any oriented link diagram $D$, the relations \eqref{eq.crossingrelations} present it as a web $\Ww$. By the above lemma, iterating the moves of Equation \eqref{eq.webrelations} will reduce $\Ww$ to a scalar multiple $\llangle \Ww\rrangle$ of the empty diagram and it is independent of the order the moves are applied. Moreover, the resulting Laurent polynomial is an invariant of the underlying oriented link $L$ and we define $\llangle L\rrangle$, the \emph{$\mathfrak{sl}_3$ invariant of $L$}, to equal $\llangle D\rrangle$ for any diagram $D$ of $L$.

\subsection{Seifert graphs, states, and the weighting}\label{sec:prelim}
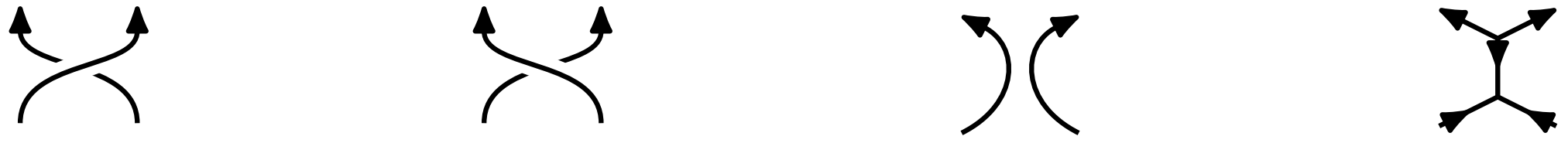
\begin{figure}[H]
\centering
\begin{subfigure}{.24\textwidth}
  \centering
    \begin{tikzpicture}[scale = 1]
            \draw[->] (1,0) to [out=90,in=-90] (0,1);
            \draw[over,->] (0,0) to[out=90,in=-90] (1,1);
    \end{tikzpicture}
  \caption*{Positive crossing}\label{fig:pcross}
\end{subfigure}%
\begin{subfigure}{.24\textwidth}
  \centering
  \begin{tikzpicture}[scale = 1]
        \draw[->] (0,0) to[out=90,in=-90] (1,1);
        \draw[over,->] (1,0) to [out=90,in=-90] (0,1);
    \end{tikzpicture}
  \caption*{Negative crossing}\label{fig:ncross}
\end{subfigure}
\begin{subfigure}{.24\textwidth}
  \centering
  \begin{tikzpicture}
			\draw[->] (0,0) .. controls (.5,.25) and (.5,.75) .. (0,1);
			\draw[->] (1,0) .. controls (.5,.25) and (.5,.75) .. (1,1);
		\end{tikzpicture}
  \caption*{{$O$-resolution}}
\end{subfigure}
\begin{subfigure}{.24\textwidth}
  \centering
  \begin{tikzpicture}
			\draw[mid>] (0,0) to (.5,.25);
			\draw[mid>] (1,0) to (.5,.25);	
			\draw[mid<] (.5,.25) to (.5,.75);
			\draw[->] (.5,.75) to (0,1);
			\draw[->] (.5,.75) to (1,1);	
		\end{tikzpicture}
  \caption*{{$W$-resolution}}
\end{subfigure}
\caption{A positive crossing, a  negative crossing, and their oriented ($O$) and web ($W$)  -resolutions.}
		\label{fig.res}
\end{figure}

With the exception of Section \ref{sec:braids}, in this paper we consider positive alternating links. An oriented link $L$ is \emph{positive} if it admits a diagram presentation $D$ with only positive crossings. Recall that the sign of a crossing is determined by the right hand rule, as seen in Figure \ref{fig.res}. 
An oriented link is \emph{alternating} if it admits an alternating diagram, i.e.\! a diagram in which any traversal along its orientation alternates between over-crossings and under-crossings.
A \emph{resolution} of a crossing is either the oriented ($O$) smoothing or the web ($W$) smoothing appearing in Figure \ref{fig.res}.
    
A \emph{state} $s$ is a choice of resolution at every crossing of $D$. We write $\alpha_+(s)$ and $\beta_+(s)$ (resp. $\alpha_-(s)$ and $\beta_-(s)$) to denote the number of $O$ and $W$ -resolutions of positive (resp. negative) crossings. We will write $\Ww(s)$ for the web diagram obtained by replacing each crossing with the indicated resolution. The scalars appearing in \eqref{eq.crossingrelations} are incorporated in the phase of $s$.

\begin{definition}
Given a state $s$, we define the phase of $s$ to be\[\phi(s):=(-1)^{\beta_+(s)+\beta_-(s)}q^{-2(\alpha_+(s)-\alpha_{-}(s))}q^{-3(\beta_+(s)-\beta_-(s))}.\] 
If $D$ is positive, we simply write $\alpha(s)$ and $\beta(s)$ and

\[
\phi(s) := (-q)^{-2\alpha(s)-3\beta(s)} = (-1)^{\beta(s)}q^{-2\alpha(s)-3\beta(s)}.
\]
\end{definition}

Due to the above discussion, we have the following proposition.

\begin{proposition}{\label{prop:statesum}}
    Let $D$ be a diagram of a link $L$. Then we can compute the $\mathfrak{sl}_3$ invariant as the sum over states \[\llangle L \rrangle = \sum_sY(s), \quad \text{where }Y(s) := \phi(s)\llangle \Ww (s)\rrangle.\]
\end{proposition}

For a state $s$ of a link diagram $D$, we have an associated \emph{state diagram} $D_s$, constructed in the following way. For $D$, consider the All-$O$ resolution, which yields disjoint circles. At each crossing in $D$ where we have a $W$-resolution, draw an edge between circles connected by a $W$-resolution. The corresponding \emph{state graph} $\Gg_s$ is obtained from $D_s$ by shrinking down each circle into a vertex, and the \emph{reduced state graph} $\Gg'_s$ is obtained by removing all multiplicity edges. The {Seifert graph} of the diagram, is the All-$W$ state graph $\Gg_W$. \\

For a graph $\Gg$, let $V(\Gg)$ and $E(\Gg)$ denote its set of vertices and edges, respectively. We will let 
\begin{align*}
    v=|V(\Gg_W)|, && e=|E(\Gg_W)|, && \mbox{and} && e'=|E(\Gg_W')|
\end{align*}
denote the number of vertices, edges, and reduced edges of $\Gg_W$. Note that every state graph is a spanning subgraph of the Seifert graph.

\section{Towards a characterization of \texorpdfstring{$\mathfrak{sl}_3$}{sl3} invariant for positive alternating links}\label{sec:posalt}
\subsection{Degree and the sorting lemma}\label{sec:sortinglemma}

Here, we state important background results from \cite{HK} that characterize the $\mathfrak{sl}_3$ invariant, and how to compute its coefficients.

\begin{definition}
 We say that $s$ and $s'$ are related by an \emph{$OW$-move} if $s'$ can be obtained from $s$ by replacing a single $O$-resolution with a $W$-resolution.
\end{definition}

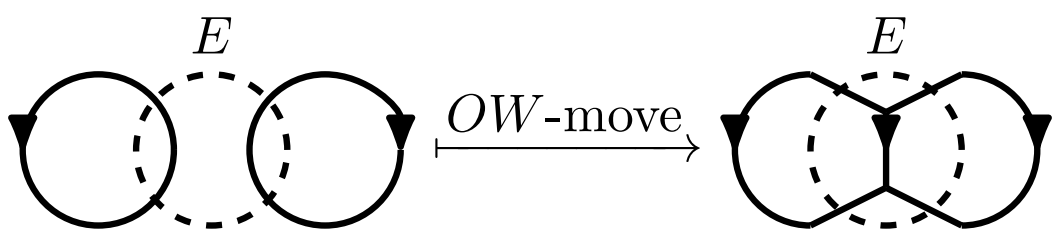
\begin{figure}[ht!]
		\[
		\begin{tikzpicture}[scale=0.5]
			\draw[mid>] (0,0) arc(0:360:1);
			\draw[->] (3,0) arc(360:0:1);
			\draw[dashed] (1.5,0) arc(0:360:1);
			\node[above] at (.5,1) {$E$};
			\node[below] at (.5,-1) {\phantom{$E$}};
		\end{tikzpicture}
		\xmapsto{\mbox{$OW$-move}}
		\begin{tikzpicture}
			\draw[dashed] (1,.5) arc(0:360:.5);
			\node[above] at (.5,1) {$E$};
			\node[below] at (.5,0) {\phantom{$E$}};
			\draw[] (0,0) to (.5,.25);
			\draw[] (1,0) to (.5,.25);	
			\draw[mid<] (.5,.25) to (.5,.75);
			\draw[] (.5,.75) to (0,1);
			\draw[] (.5,.75) to (1,1);		
			\draw[mid<] (0,0) arc(270:90:.5);
			\draw[mid<] (1,0) arc(-90:90:.5);
		\end{tikzpicture}
		\]\vspace{-1\baselineskip}
		\caption{An $OW$-move between two state webs  $\mathbb{W}$(left) and $\mathbb{W}'$(right)  in the disc $E$ whose boundary is indicated by the dotted line.}\label{fig.OW}
			\end{figure}

A simple example of an $OW$-move is given in Figure \ref{fig.OW}.
	
For any Laurent polynomial $f(q)$, we write $\deg(f)$ to mean the degree of its leading term.

\begin{lemma}[{\cite[Prop~2.6]{HK}}]\label{lem:owstep}
Let $s$ and $s'$ be states of a positive diagram which are related by an $OW$-move. Then, $\deg(Y(s')) = \deg (Y(s))$ or $\deg(Y(s')) = \deg (Y(s))-2$.
\end{lemma}

The above lemma may also be formulated in the language of state graphs.

\begin{lemma}[{\cite[Rem~3.3]{HK}}]\label{lem:components}
    For a state $s$ of a positive diagram whose state graph $\Gg_s$ has exactly $v-k$ components, there is a degree bound $\deg(Y(s)) \le 2(v-e-k)$.
\end{lemma}

  If we consider the All-$O$ state which has $v$ components, as $\Gg_O$ has no edges, this implies that the $\inv$ invariant of a positive link $L$ has highest degree $2(v-e)$ where $v$ and $e$ are determined from any positive diagram. Since we know the degree either decreases by two or zero by introducing more $W$-resolutions, the $\mathfrak{sl}_3$ polynomial has the form
\[\llangle L \rrangle = \gamma_1q^{2(v-e)} + \gamma_2q^{2(v-e-1)} + \gamma_3q^{2(v-e-2)} + \gamma_4q^{2(v-e-3)} + \gamma_5q^{2(v-e-4)} + \text{lower degree terms.}\]
Expressions for $\gamma_1,\gamma_2,\gamma_3$ were given in terms of diagrammatic quantities of any positive diagram for $L$ in \cite{HK}.

\subsection{The fourth and fifth coefficients}\label{subsec:45}

In this section, we compute the fourth and fifth coefficients, $\gamma_4$ and $\gamma_5$, in terms of data read from the Seifert graph. First, we will introduce key lemmas needed for the computation before determining the coefficients in Subsections \ref{sec:fourth} and \ref{sec:fifth}. For now, let us assume the underlying Seifert graphs are acyclic. We will handle cycles appearing as a special case when computing the fourth and fifth coefficients explicitly. 

\begin{lemma}\label{lem:ys}
    Let $D$ be a positive link diagram and let $s$ be a state such that $\Gg_s'$ is acyclic. Then 
    \[Y(s) = (-1)^{\mathbf{k}}q^{-2e-\mathbf{k}}[3]^{\,v-r}[2]^{\mathbf{k}}\]
    where $r=|E(\Gg_s')|$, $\mathbf{k} = \sum_{i \in E(\Gg'_s)}k_i$, and where $k_i$ is the multiplicity of the reduced edge $i$ in $\Gg_s$.
\end{lemma}
\begin{proof}
    Since $D$ is positive, $\phi(s) =(-1)^{\beta}q^{-2\alpha-3\beta}$. By definition, $\beta = \mathbf{k}$ and $\alpha = e - \mathbf{k}$, thus \[\phi(s)= (-1)^{\mathbf{k}}q^{-2e-\mathbf{k}}.\]

    Next consider the web evaluation $\llangle \Ww(s) \rrangle$. Since $\Gg_s'$ is acyclic with $v$ vertices and $r$ edges, it has $v-r$ components. Consider a component of $\Gg_s$ and its associated connected web $W_1$. Every multiplicity edge and edge in $\Gg_s$ corresponds to a bubble $W_1$, which contributes a factor of $[2]$ to $\llangle W_1 \rrangle$. Since $W_1$ is connected and remains connected after applying the bubble move, this component finally contributes a factor of $[3]$ by applying a circle move. Applying this argument across all $v-r$ components of $\Ww(s)$ yields $\llangle \Ww(s) \rrangle = [3]^{v-r}[2]^{\mathbf{k}}$ and $Y(s) = (-1)^{\mathbf{k}}q^{-2e-\mathbf{k}}[3]^{\,v-r}[2]^{\mathbf{k}}$.
\end{proof}

\begin{remark}\label{rem:forest}
    Lemma \ref{lem:ys} shows that among state graphs of a positive link diagram $L$ which are acyclic, their web evaluations depend only on the number of components and edges of the reduced graph. For example,
        \[
    \leftlangle\,\,\,
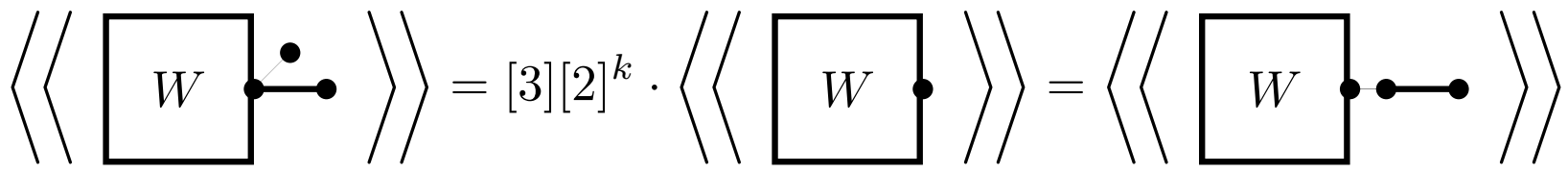
    \begin{tikzpicture}
        \node (W)[draw, minimum size=1cm] {$W$};
        \fill[black] (W.east) circle (2pt);
        \draw (W.east) -- ++(0.5,0);
        \fill[black] (W.east) -- ++(0.5,0) circle (2pt);
        \fill[black] (W.east) -- ++(0.25,0.25) circle (2pt);
    \end{tikzpicture}
    \,\,\,\rightrangle=[3][2]^{k}\cdot
    \leftlangle\,\,\,
    \begin{tikzpicture}
        \node (W)[draw, minimum size=1cm] {$W$};
        \fill[black] (W.east) circle (2pt);
    \end{tikzpicture}
    \,\,\,\rightrangle=
    \leftlangle\,\,\,
    \begin{tikzpicture}
        \node (W)[draw, minimum size=1cm] {$W$};
        \fill[black] (W.east) circle (2pt);
        \fill[black] (W.east) -- ++(0.25,0) circle (2pt);
        \fill[black] (W.east) -- ++(0.75,0) circle (2pt);
        \draw (0.75,0) to (1.25,0);
    \end{tikzpicture}
    \,\,\,\rightrangle
    \]

    where $k$ is the multiplicity of the added edge.
\end{remark}

\begin{lemma}\label{lem:cont} Let $s$ be a state where the reduced state graph $\Gg'_s$ is a forest. Let $r$ be the number of edges in $\Gg'_s$, and $\mathbf{k} = k_1 + \cdots +k_{r}$ where $k_i$ is the multiplicity of edge $i$ in $\Gg_s$. The coefficient contribution to the $2(v-e-n)$ degree term of the $\mathfrak{sl}_3$ invariant from such a state is given by \[\sum_{i+j+m=n-r}(-1)^{\mathbf{k}}\binom{v-r}{i}{\binom{i}{j}}\binom{\mathbf{k}}{m}.\]
\end{lemma}

\begin{proof}

We know that for a state $s$ as described above, we have \[Y(s) = (-1)^{\mathbf{k}}q^{-2e- \mathbf{k}}[3]^{v-r}[2]^{\mathbf{k}}.\]

Expanding $[3]^{v-r}$ and $[2]^{\mathbf{k}}$, we get 

\[Y(s) = (-1)^{\mathbf{k}}q^{-2e- \mathbf{k}}\left[\sum_{i=0}^{v-r}\binom{v-r}{i}q^{2v-2r-2i}\sum_{j=0}^i\binom{i}{j}q^{-2j}\right]\left[ \sum_{m=0}^{\mathbf{k}}\binom{\mathbf{k}}{m} q^{\mathbf{k}-2m} \right].\]

Looking for the coefficient of $q^{2(v-e-n)}$ in this expansion gives us \begin{align*}
    2(v-e-n) &= 2v-2e-2r-2i-2j-2m \iff n-r=i+j+m, \\&\qquad \text{ where } \quad 0\le j\le i \le v-r\quad \text{ and } \quad  0\le m \le \mathbf{k}.
\end{align*}
Thus, the coefficient of $q^{2(v-e-n)}$ is exactly the sum over all satisfactory tuples $(i,j,m)$.
\end{proof}

\begin{lemma}\label{lem:delta}
For $n \ge 1$ and $m \ge 0$,
\[
\sum_{k=1}^{n}\binom{n}{k}(-1)^{k}\binom{k}{m} = (-1)^{m}\delta_{n,m} - \delta_{m,0}.
\]
\end{lemma}

\begin{proof}
    Let $m \in \mathbb Z_{\ge0}$. First we consider the case        $m=0$. We have \[\sum_{k=1}^n\binom{n}{k}(-1)^{k}\binom{k } {0} =-1.\] For $m > 0$, by the binomial formula, \[(x+1)^n = \sum_{k=0}^n\binom{n}{k}x^k.\] Taking the $m$-th derivative of both sides, we get \[n(n-1)(n-2)\cdots(n-m+1)(x+1)^{n-m} = \sum_{k=0}^n\binom{n}{k}x^{k-m}k(k-1)\cdots(k-m+1).\]

    Substituting $x=-1$ yields \[n(n-1)(n-2)\cdots(n-m+1)(0)^{n-m} = \sum_{k=0}^n\binom{n}{k}(-1)^{k-m}k(k-1)\cdots(k-m+1).\]

    Note that for $0\le k\le m-1$, the terms of the summations vanish. To match the summation in our theorem statement, we multiply both sides by $(-1)^m/m!$, and since the $k=0$ term is $0$, we get

    \[(-1)^m\cdot\frac{n(n-1)(n-2)\cdots(n-m+1)}{m!}\cdot(0)^{n-m} = 0+\frac{1}{m!}\sum_{k=1}^n\binom{n}{k}(-1)^{k}k(k-1)\cdots(k-m+1) \] \[= \sum_{k=1}^n\binom{n}{k}(-1)^{k}\binom{k}{m}.\]

    Now consider \[(-1)^m\cdot\frac{n(n-1)(n-2)\cdots(n-m+1)}{m!}\cdot(0)^{n-m}.\] We have two cases:
    \begin{itemize}
        \item[(a)] For $n \ne m$, this expression is $0$ as $0^{n-m} = 0$ for $n \ne m$.
        \item[(b)] For $n= m$, this expression reduces to \[(-1)^m\cdot\frac{m(m-1)(m-2)\cdots(1)}{m!}\cdot(1) = (-1)^m.\]
    \end{itemize}

    So, if $m=0$, the expression equals $-1$. For $m >0$, we have $(-1)^m\delta_{n,m}$. This proves the lemma for all $m \in \mathbb Z_{\ge0}$.
\end{proof}

\subsubsection{The fourth coefficient}\label{sec:fourth}

We organize states by the number of connected components that they have relative to the All-$O$ state, based on Lemma \ref{lem:components}. For states whose state graphs are acyclic, we may apply Lemma \ref{lem:cont}. We will consider state graphs which contain a cycle separately. In which case, we may apply Remark~\ref{rem:forest} to describe the contribution of a number of states.\\

First, we will look into the case of acyclic state graphs. From Lemma \ref{lem:cont}, we know the coefficient contribution from a single state $s$. In order to find the total contribution from all the states to a given degree, we first need to iterate over the different multiplicities of each edge in the reduced graph, as two states with different multiplicities between two vertices in $\Gg_s$ have the same reduced graphs. Hence, we must account for all the different combinations of the multiplicities for each edge in $\Gg_s'$.\\ 

Say there are $r$ edges in the reduced graph of a state, and let $z_i$ denote the multiplicity of edge $i$ for $1\leq i\leq r$. Since each state graph is a spanning subgraph of the Seifert graph, we need to account for all the different combinations of edges in $\Gg_s'$ that have $r$ edges. The following lemma formalizes this into an expression.

\begin{lemma}\label{lem:counting}
    Let $D$ be a positive link diagram. The contribution of states with acyclic reduced state graphs to the $2(v-e-n)$ degree term of the $\mathfrak{sl}_3$ invariant is given by 
    \[
    \sum_{r=0}^{e'}\;
\sum_{\substack{S \subseteq E(\Gg'_W) \\ |S| = r\\\Gg'_s \text{is acyclic}}}\;
\sum_{k_{i_1}=1}^{z_{i_1}}\cdots \sum_{k_{i_r}=1}^{z_{i_r}}\left(\prod_{\ell \in S}\binom{z_\ell}{k_\ell}\right)
\sum_{i+j+m=n-r}(-1)^{\mathbf{k}}\binom{v-r}{i}\binom{i}{j}\binom{\mathbf{k}}{m}
    \]
    where $z_i$ is the multiplicity of the $i$-th reduced edge of $\Gg_W$.
\end{lemma}
\begin{proof}
    From Lemma \ref{lem:cont}, each state $s$ with acyclic state graph contributes 
    \[
    \sum_{i+j+m=n-r}(-1)^{\mathbf{k}}\binom{v-r}{i}{\binom{i}{j}}\binom{\mathbf{k}}{m}
    \]
    to the coefficient of degree $2(v-e-n)$, where $r=|E(\Gg'_s)|$, $k_i$ is the multiplicity of edge $i$ in $\Gg_s$, and $\mathbf{k} = k_1 + \cdots + k_{r}$. \\
    
    It remains to enumerate over all states. Fix an edge support set $S=E(\Gg_s')\subseteq E(\Gg_W')$. Now there are $\prod_{\ell\in S}\binom{z_\ell }{k_\ell}$ states with the same underlying reduced state graph and support $S$, considered as distinct subgraphs of $\Gg_W$. Each of these have equal contribution as the state considered in the previous paragraph. This yields 
    \[
    \left(\prod_{\ell \in S}^{z_\ell}\binom{z_\ell}{k_\ell}\right)
\sum_{i+j+m=n-r}(-1)^{\mathbf{k}}\binom{v-r}{i}\binom{i}{j}\binom{\mathbf{k}}{m}.
\]
    Next, we vary the multiplicity of these edges by allowing each $k_\ell$ to vary from $1$ to $z_\ell$. Finally summing over all support sets $S$ of size $r$, then allowing $r$ to vary between $0$ and $e'=E(\Gg_W')$ produces the desired sum. 
\end{proof}

Finding each integral solution to $i+j+m=n-r$, we can explicitly evaluate the innermost sum, and since these variables are independent of the edges chosen, the sum can be factored out. We can expand $\binom{\mathbf{k}}{m}$, by applying Vandermonde's Identity iteratively. Recall that this identity says for all nonnegative integers $m,n$ and $r$, the following holds:
    \[\sum_{k=0}^{r} \binom{m}{k} \binom{n}{r-k} = \binom{m+n}{r}.\]
Now, match the resulting terms to the form of Lemma \ref{lem:delta} to make the entire summation independent of $S$, turning the sum over subsets into multiplying by a factor of $\binom{e'}{r}$. Once, we are in this form, applying Lemma $\ref{lem:delta}$ yields Kronecker deltas $\delta_{z_\ell,m}$. As this value depends only on the multiplicity of the edge $\ell$, we introduce the statistic
\begin{equation}\label{eq:edge_count}
    \mu_{m}=|\{\ell \in E(\mathbb{G}'_{W}) : z_\ell \neq m\}|.
\end{equation}

So far, we have discussed forests specifically because the existence of cycles in the reduced state graph cause the web evaluations to not be easily generalized. The first appearance of a cycle in a state graph occurs when considering states with $v-3$ components. A direct computation of such a state, call it $s_c$, gives 
\begin{equation}\label{eq:sc}
    Y(s_c)=2(-1)^{k_1+k_2+k_3+k_4}q^{-2e-k_1-k_2-k_3-k_4}[3]^{(v-3)}[2]^{k_1+k_2+k_3+k_4-2}
\end{equation}
where $k_i$ is the multiplicity of edge $i$ in $\Gg_{s_c}$. However, we observe that the highest degree term of this polynomial is in fact $2(v-e-4)$, and therefore the four-cycle configuration does not contribute to the fourth coefficient despite having $p-3$ components. Instead, it contributes to $\gamma_5$,  see Proposition \ref{prop:g5-cycle}.

\begin{proposition}\label{prop:g4} The following states contribute to $\gamma_4$:
    \begin{itemize}
        \item the All-$O$ state: $2\binom{v}{2}+\binom{v}{3}$,
        \item the $v-1$ component states: $-e'\binom{v}{2}+(1-v)(e'-\mu_1)+(e'-\mu_2)$,
        \item the $v-2$ component states: $\binom{e'}{2}(v-2) + (e'-\mu_1)\mu_1 + 2\binom{e'-\mu_1}{2}$,
        \item the $v-3$ component states: $-\binom{e'}{3}$ to $\gamma_4$. 
    \end{itemize}
\end{proposition}

The proof of the above proposition follows by considering fixed $r$ sums in Lemma \ref{lem:counting}, where $r$ varies between $0$ and $3$. The sum over these terms gives the total contribution to the $2(v-e-3)$ term, hence a formula for $\gamma_4$. 

\begin{theorem}\label{thm:gamma4}
For any positive alternating link diagram,
\[
\gamma_4 = \binom{v-e'+2}{3} + (\mu_1-1)\,(v-e') + e' - \mu_2.
\]

\end{theorem}

\subsubsection{The fifth coefficient}\label{sec:fifth}
The same methods describe above, now considering states with $v$ through $v-4$ components of $\Gg'_s$, nearly computes $\gamma_5$. However, we must also account for four-cycles in state graphs. We let $\sigma$ denote the number of four-cycles in $\Gg_W'$. We organize the states by the number of components of $\Gg'_s$.
\begin{proposition}\label{prop:g5}
The following states contribute to $\gamma_5$:
\begin{itemize}
    \item the All-$O$ state: $\binom{v}{2}+v\binom{v-1}{2}+\binom{v}{4}$,
    \item the $v-1$ component states: {\footnotesize$-(e'-\mu_3)+(v-1)(e'-\mu_2)-\Big((v-1)+\tbinom{v-1}{2}\Big)(e'-\mu_1)-\Big((v-1)(v-2)+\tbinom{v-1}{3}\Big)e'$},
    \item the $v-2$ component states: {\footnotesize\[\binom{e'-\mu_1}{2}-(e'-\mu_2)\mu_2-2\binom{e'-\mu_2}{2}-2(e'-\mu_1)\mu_1-4\binom{e'-\mu_1}{2}+v(e'-\mu_1)\mu_1+2v\binom{e'-\mu_1}{2}+\Big(v-2+\tbinom{v-2}{2}\Big)\binom{e'}{2},\]}
    \item the $v-3$ component states: {\footnotesize$-\binom{e'}{3}(v-3)-(e'-\mu_1)\binom{\mu_1}{2}-2\binom{e'-\mu_1}{2}\mu_1-3\binom{e'-\mu_1}{3}$},
    \item the $v-4$ component acyclic states: $\binom{e'}{4} - \sigma$.
\end{itemize}
\end{proposition}

Again, these contributions were calculated using Lemma \ref{lem:counting}. In the last item above, we subtract $\sigma$, as cycles do not meet the hypothesis of the lemma. Each choice of four edges in the Seifert graph contributes $1$ to the $\gamma_5$ except when a cycle is generated. We handle that contribution separately.
\begin{proposition}\label{prop:g5-cycle}
The $v-3$ component states containing a four-cycle contribute $2\sigma$ to $\gamma_5$.
\end{proposition}

\begin{proof}
    Let $s_c$ be some state containing a four-cycle. Equation \eqref{eq:sc} gives an expression for $Y(s_c)$ with leading term $2(-1)^{k_1+k_2+k_3+k_4}q^{2(v-e-4)}$. Summing over edge multiplicities $k_i$ varying from $1$ to $z_i$, the multiplicity of edge $i$ in $\Gg_W$, we have
    \[
    2\sum_{k_1=1}^{z_1}\sum_{k_2=1}^{z_2}\sum_{k_3=1}^{z_3}\sum_{k_4=1}^{z_4}(-1)^{k_1+k_2+k_3+k_4}\binom{z_1}{k_1}\binom{z_2}{k_2}\binom{z_3}{k_3}\binom{z_4}{k_4}=2.
    \]
Then, we take this over all four-cycles, giving $2\sigma$.
\end{proof}

Now summing the contributions of Proposition \ref{prop:g5} and \ref{prop:g5-cycle}, we have a formula for the fifth coefficient of the invariant, which can be simplified to the form shown below. 

\begin{theorem}\label{thm:gamma5}
For any positive alternating link diagram,
\[
\gamma_5 = \binom{v-e'+3}{4} + (\mu_1-2)\binom{v-e'+1}{2} + (e'-\mu_2+1)(v-e') + \binom{\mu_1}{2} - e' + \mu_3 + \sigma 
\,.
\]
\end{theorem}

\subsection{When \texorpdfstring{$\Gg'_W$}{the reduced Seifert graph} is a tree}

In this section, we present a characterization of the invariant for a link $L$ which admits a positive alternating diagram in the case $\Gg_W'$ is a tree. Equivalently, by \cite[Theorem 5.11]{Gutsbook}, $L$ is a positive alternating fibered link. Moreover the class of fibered positive alternating links are the links which can be written as connected sums of $(2,n)$ torus links \cite[Corollary 7.2]{HK}.
 
\begin{theorem}\label{thm:fibered}
    For every link $L$ with a positive alternating connected link diagram, whose reduced Seifert graph $\Gg_W'$ is a tree, we have
\begin{equation}{\label{eqn:counting}}
    \llangle L \rrangle = q^{-2e}\sum_{r=0}^{e'}[3]^{\,v-r}\,e_r(u_1,\dots,u_{e'}),
  \qquad u_i = (-q^{-2})^{z_i}-1,
\end{equation}
where $e_r$ is the $r$-th elementary symmetric polynomial and $z_i$ is the multiplicity of the $i$-th reduced edge.
\end{theorem}

\begin{proof}
    Let $S$ be a subset of $E(\Gg'_W)$. Each state $s$ satisfies $Y(s) = (-1)^{\mathbf{k}}q^{-2e-\mathbf{k}}[3]^{\,v-r}[2]^{\mathbf{k}}$ by Lemma \ref{lem:ys} with $\mathbf{k} = \sum_{\ell \in S} k_\ell$.  We can rewrite $Y(s)$ equivalently as \[Y(s) = q^{-2e}[3]^{v-r}\prod_{\ell \in S}\left((-1)^{k_\ell}q^{-k_{\ell}}[2]^{k_\ell}\right), \qquad r=|S|.\] Furthermore, since $(-1)^{k_\ell}q^{-k_\ell}[2]^{k_\ell}=(-1-q^{-2})^{k_\ell}$, we have \[Y(s) = q^{-2e}[3]^{v-r}\prod_{\ell \in S}(-1-q^{-2})^{k_\ell}.\]

    Now, to sum over all possible states, similar to Lemma $\ref{lem:counting}$, we have 
    \[\sum_{r=0}^{e'}\quad\sum_{\substack{S\subseteq E(\Gg'_W)\\|S|=r}} \left(\sum_{k_{i_1}=1}^{z_{i_1}}\cdots\sum_{k_{i_r}=1}^{z_{i_r}}\binom{z_{i_1}}{k_{i_1}}\cdots \binom{z_{i_r}}{k_{i_r}} q^{-2e}[3]^{v-r}\prod_{\ell \in S}(-1-q^{-2})^{k_\ell}\right).\]

    We can continue to simply as \[q^{-2e}\sum_{r=0}^{e'}[3]^{v-r}\quad\sum_{\substack{S\subseteq E(\Gg'_W)\\|S|=r}}\prod_{\ell\in S}\left(\sum_{k_\ell = 1}^{z_\ell}\binom{z_\ell}{k_\ell}(-1-q^{-2})^{k_\ell}\right).\]

    Applying the binomial theorem gives us
    \[\prod_{\ell \in S}\left[\sum_{k_\ell=1}^{z_\ell}\binom{z_\ell}{k_\ell}(-1-q^{-2})^{k_\ell}\right]=\prod_{\ell \in S}\left((-q^{-2})^{z_\ell}-1\right).\]
    Thus, we have

    \[q^{-2e}\sum_{r=0}^{e'}[3]^{v-r} \sum_{\substack{S \subseteq E(\Gg'_W) \\ |S| = r}}\prod_{\ell \in S}\left((-q^{-2})^{z_\ell}-1\right).\]
    
    Substituting $u_i := (-q^{-2})^{z_i}-1$ and identifying the latter sum as a symmetric function \cite[Chapter 9, Equation 9.2]{martin} gives the result.
\end{proof}

\begin{example}
    We will compute the invariant for the trefoil using Theorem \ref{thm:fibered}. The Seifert graph has two vertices and one reduced edge with multiplicity three. This gives 
    \[q^{-6}\sum_{r=0}^{1}[3]^{2-r}e_r(x_1) = q^{-6}\big([3]^{2}+[3]\cdot u_1\big) = q^{-6}\big([3]^{2}+[3](-q^{-6}-1)\big) = q^{-2}+q^{-4}+2q^{-6}+q^{-8}-q^{-12}-q^{-14}\]
\end{example}

As an application of Theorem \ref{thm:fibered}, we can express the invariant over the same class of links as a Tutte polynomial evaluation. For a reference on matroids and the Tutte polynomial, one can refer to the textbook \cite[Chapters 3-4]{martin}. We will let $T(M;x,y)$ denote the Tutte polynomial of the matroid $M$ evaluated at $x$ and $y$. 

\begin{theorem}
    Let $L$ be a connected link which admits a positive alternating diagram, and whose reduced Seifert graph $\Gg_W'$ is a tree. Then \[\llangle L\rrangle = q^{-2e}[3](-q^{-2}-1)^{e'}T(M(\Gg_W); \frac{-q^4}{q^{2}+1}, -q^{-2}).\]
\end{theorem}

\begin{proof}
    Since $\Gg_W'$ is a tree, $\Gg_W$ is a connected graph of parallel classes. Thus, \[M(\Gg_W) \cong \bigoplus_{\ell\in E(\Gg_W')} U_1(z_\ell)\]

    where $U_k(n)$ is the uniform matroid of rank $k$ on $n$ elements. Since the Tutte polynomial is multiplicative over direct sums, we have \[T(M(\Gg_W);x,y) = \prod_{\ell\in E(\Gg_W')} T(U_1(z_\ell);x,y).\]

    The Tutte polynomial of $U_1(z_\ell)$ is known:
    \[T(U_{1}(z_\ell);x,y) = x+y+y^2+\dots+y^{z_{\ell}-1}.\]

    Thus, \[T(M(\Gg_W);x,y) = \prod_{\ell\in E(\Gg_W')} (x+y+y^2+\dots+y^{z_\ell-1}).\]

    We further simplify as \[\prod_{\ell \in E(\Gg_W') }\left(x+\frac{y^{z_\ell}-y}{y-1}\right) = (y-1)^{-e'}\prod_{\ell \in E(\Gg'_W)}(x(y-1)+y^{z_\ell}-y) = (y-1)^{-e'}\prod_{\ell \in E(\Gg_W')}((x-1)(y-1)+y^{z_\ell}-1)\]

    Now, consider the Equation $\eqref{eqn:counting}$. We can rewrite the equation into a product using the generating function for elementary symmetric polynomials: \[\sum_{k=0}^nt^{k}e_k(x) = \prod_{i=1}^n(1+tx_i)\] which can be seen in \cite[Chapter 9, Equation 9.2]{martin}. Using this, we can rewrite Equation $\eqref{eqn:counting}$ as 
    \begin{align*}
        &q^{-2e}\sum_{r=0}^{e'}[3]^{\,v-r}\,e_r(u_1,\dots,u_{e'}) = q^{-2e}[3]^{e'+1}\sum_{r=0}^{e'}[3]^{-r}e_r(u_1,\dots,u_{e'})=q^{-2e}[3]^{e'+1}\prod_{i=1}^{e'}(1+[3]^{-1}u_i) \\
        &=q^{-2e}[3]^{e'+1}\prod_{\ell \in E(\Gg'_w)}(1+[3]^{-1}u_\ell)=q^{-2e}[3]\prod_{\ell \in E(\Gg_W')}([3]+u_\ell)=q^{-2e}[3]\prod_{\ell \in E(\Gg_W')}([3]+(-q^{-2})^{z_\ell}-1).
    \end{align*}
    
    We want to match this product over reduced edges to the expression of the Tutte polynomial we found earlier. From the similarity in the product expression for the Tutte polynomial and the invariant, the natural substitution is $y=-q^{-2}$, and $x$ is then forced:

    \[(x-1)(y-1) = [3] \iff x=[3](y-1)^{-1}+1 = \frac{q^2}{-q^{-2}-1}\]

    Accounting for the factor $(y-1)^{-e'}$ taken out of the Tutte polynomial, we find the proposed equality.
\end{proof}

\section{Behavior under connected sums}\label{sec:connect}
In this section, we show that the fourth and fifth coefficients behave nicely under connected sums. Throughout, let $e'(D_i)$ denote the number of edges in the reduced Seifert graph of the corresponding knot $K_i$, and similarly $\sigma(D_i)$ denote the number of four-cycles.

\begin{theorem}[{\cite[Theorem 6.1]{HK}}]\label{thm:HKp} Let $K$ be a positive knot with prime factors $K_1, \dots, K_p$.
\begin{itemize}
    \item[(1)] $\gamma_2(K) = 1 - p + \sum_{j=1}^p\gamma_2(K_j)$
    \item[(2)] $\lambda(K) = \sum_{j=1}^p\lambda(K_j), \quad\text{where}\quad \lambda(K_j) := \mu(K_j)-\theta(K_j)$.
    \item[(3)] $\gamma_3(K) = \frac{(\gamma_2(K)+1)\gamma_2(K)}{2} + \lambda_1(K)$ 
\end{itemize}
\end{theorem}

Note that the conventions in the above theorem differ from ours. In the work of \cite{HK}, their $\mu$ is equivalent to our definition of $\mu_1$. We have not defined $\theta$, which is the number of mixed states, since in the alternating case, $\theta = 0$. So, we define $\lam_i := \sum_{j=1}^p\mu_i(K_j)$.
Due to {\cite[Theorem 1]{Men}}, we know a non-prime alternating link $L$ factors as a connected sum of non-trivial alternating links. We can iterate this until we reach prime factors.\\

Define $\epsilon := \sum_{j=1}^p e'(D_j)$ and $\varsigma := \sum_{j=1}^p \sigma(D_j)$. Using the above definition of $\gamma_2(K)$ in Theorem \ref{thm:HKp}(1), we have the following result.

\begin{theorem}\label{thm:connectsum} Let $K$ be a positive alternating knot with prime factors $K_1, \dots, K_p$.
\begin{itemize}
    \item[(1)] \label{itm:cntsum1}$\gamma_4(K)=\binom{\gamma_2(K)+2}{3} + (\lam_1-1)\,(\gamma_2(K)) + \epsilon - \lam_2$
    \item[(2)] \label{itm:cntsum2}$\gamma_5(K) =\binom{\gamma_2(K)+3}{4} + (\lam_1-2)\binom{\gamma_2(K)+1}{2} + (e'-\lam_2+1)\gamma_2(K) + \binom{\lam_1}{2} - \epsilon + \lam_3 + \varsigma$
\end{itemize}
\end{theorem}

\begin{proof}
    We proceed by induction on $p$. For the base case $p=1$, we have our original formulae in Theorems \ref{thm:gamma4} and $\ref{thm:gamma5}$. Assume that the claim holds for all positive alternating knots with prime decomposition number at most $p$. Let $K = \#_{j=1}^pK_j$ be such a knot, and let $\Kp$ be a prime positive alternating knot. Let $D$ and $D_{p+1}$ be positive alternating diagrams of $K$ and $\Kp$, with vertex and edge counts $v,e'$ and $v_{p+1}, e_{p+1}'$ respectively. We will write $\gamma_i(K_j)$ to denote the $i$th coefficient of $\llangle K_j\rrangle$.
    The HOMFLY skein relation specialized to the $\mathfrak{sl}_3$ invariant, gives \[[3]\llangle K \# K_{p+1} \rrangle = \llangle K \rrangle \llangle K_{p+1} \rrangle.\] 
    First, we will compute the right hand side. These factors have the form
    \begin{align*}
        &\llangle K \rrangle=q^{2(v-e')}+\gamma_2(K)q^{2(v-e'-1)}+\gamma_3(K)q^{2(v-e'-2)} + \text{ lower degree terms}\\
        &\llangle K_{p+1} \rrangle=q^{2(v_{p+1}-e_{p+1}')}+\gamma_2(\Kp)q^{2(v_{p+1}-e'_{p+1}-1)}+\gamma_3(\Kp)q^{2(v_{p+1}-e'_{p+1}-2)} + \text{ lower degree terms}
    \end{align*}
Expanding the product gives us
    \[\llangle K \rrangle \llangle K_{p+1} \rrangle=q^{2(v'' - \,e'')}+\gamma_2'q^{2(v''- \,e'' - \,1)}+\gamma_3'q^{2(v''- \,e''- \,2)}+ \text{ lower order terms}\]
where $v'' := v+v_{p+1},\, e'' := e'+e'_{p+1}$.
From this expansion, we have 
\begin{align*}
    &\gamma_2' = \gamma_2(K)+\gamma_2(\Kp)\\
    &\gamma_3' = \gamma_3(K)+\gamma_3(\Kp)+ \gamma_2(K)\gamma_2(\Kp)\\
    &\gamma_4'=\gamma_4(K)+\gamma_4(K_{p+1})+\gamma_3(K)\gamma_2(\Kp) + \gamma_2(K)\gamma_3(\Kp)\\
    &\gamma_5'=\gamma_5(K)+\gamma_5(\Kp)+\gamma_4(K)\gamma_2(\Kp)+\gamma_2(K)\gamma_4(\Kp)+\gamma_3(K)\gamma_3(\Kp)
\end{align*}
Expanding $[3] = q^2+1+q^{-2}$ and comparing coefficients on both sides gives 
    \begin{align*}
        &\gamma_2(K \# K_{p+1})=\gamma_2'-1\\
        &\gamma_3(K \# K_{p+1})=\gamma_3'-\gamma_2'\\
        &\gamma_4(K \# K_{p+1})=\gamma_4'-\gamma_3'+1\\
        &\gamma_5(K \# K_{p+1})=\gamma_5'-\gamma_4'+\gamma_2'-1
    \end{align*}
Now substitute $\gamma_5',\gamma_4',$ and $\gamma_3'$ into each equation. Let $\lam_i':= \lambda_i(K) + \lambda_i(\Kp)$, and $\sigma' = \sigma(D) + \sigma(D_{p+1})$. Collecting binomial coefficients using Vandermonde's identity and Pascal's identity gives
\begin{align*}
\gamma_3' + 1 &= \binom{\gamma_2' + 1}{2} + \lam_1' + 1,\\
\gamma_4' &= \binom{\gamma_2' + 1}{3} + \binom{\gamma_2' + 1}{2} - \gamma_2'
    - \lam_1' + (\gamma_2' + 1)\lam_1' - \lam_2' + e'',\\
\gamma_5' &= \binom{\gamma_2' + 1}{4} + 2\binom{\gamma_2' + 1}{3}
    - \binom{\gamma_2'}{2} + \binom{\gamma_2' + 1}{2}\lam_1'
    + \binom{\lam_1'}{2} - \gamma_2'\lam_2' + (\gamma_2' - 1)\,e''
    + \lam_3' + \sigma'.
\end{align*}
Using Pascal's identity further gives
\begin{align*}
\gamma_4(K \# \Kp)
&= \binom{\gamma_2'}{3} + \binom{\gamma_2'}{2}  - (\gamma_2' - 1)- \lam_1'
    + \gamma_2'\lam_1' - \lam_2' + e''\\
    &= \binom{\gamma_2' + 1}{3} - (\gamma_2' - 1) - 2\lam_1'
    + (\gamma_2' + 1)\lam_1' - \lam_2' + e''
\end{align*}
Note that $\lambda_i' = \lambda_i(K \# \Kp)$, and $e'' = \epsilon(K \# \Kp)$, which completes the proof of Theorem \ref{itm:cntsum2}(1). Similarly, for the next coefficient, we have
\begin{align*}
    \gamma_5(K\#\Kp) &=  \binom{\gamma_2'}{4} + 2\binom{\gamma_2'}{3}
    - \binom{\gamma_2' - 1}{2} + \binom{\gamma_2'}{2}\lam_1'
    + \binom{\lam_1'}{2} - (\gamma_2' - 1)\lam_2' + (\gamma_2' - 2)\,e''
    + \lam_3' + \sigma'\\
    &= \binom{\gamma_2' + 2}{4}
    + (\lam_1' - 2)\binom{\gamma_2'}{2}
    + (e'' - \lam_2' + 1)\,(\gamma_2'-1) + \binom{\lam_1'}{2} - e''
    + \lam_3' + \sigma'
\end{align*}
Again, note that $\sigma' = \varsigma(K \# \Kp)$, thus we have shown Theorem \ref{itm:cntsum2}(2), which completes the induction.
\end{proof}
We now use these results to analyze the class of knots that can be represented as the closure of positive alternating braids.

\begin{theorem}[{\cite[Theorem 7.1]{HK}}]\label{thm:harper}
If a knot $K$ of prime decomposition number $p$ is a positive closed braid, then $\gamma_1(K)=\gamma_2(K)=1$ and $\gamma_3(K)=p+1$.
\end{theorem}
It was shown in \cite[Corollary 7.2]{HK} that the only nontrivial positive alternating braid closures are connected sums of $(2,n)$ torus links. We can express our formulas in terms of higher degree coefficients, and directly extend this result to ours.
\begin{theorem}{\label{thm:toruslink}}
    Let $K$ be a knot with prime decomposition number $p$ which is an alternating positive closed braid, and thus is a connected sum of $(2,n)$ torus knots. Then
    \begin{align*}
        &\gamma_4(K) = p + \epsilon- \lambda_2\\
        &\gamma_5(K) = \binom{p+1}{2} +\lam_3-\lam_2+\varsigma
    \end{align*}
\end{theorem}
\begin{proof} Let $K$ be a knot with prime decomposition number $p$, with positive alternating diagram $D$. First, we will present our coefficient formulae in an alternative form in terms of the coefficient of one lower index, a similar form as given in the proof of Theorem $\ref{thm:connectsum}$.
\begin{align*}
    &\gamma_4(K) = \binom{\gamma_2(K)+1}{3}+\gamma_3(K)-\gamma_2(K)+(\gamma_2(K)-1)\lam_1-\lam_2+\epsilon\\
    &\gamma_5(K) = \binom{\gamma_2(K)+2}{4}+(\lam_1 - 2)\binom{\gamma_2(K)}{2} + \binom{\lambda_1}{2} + \gamma_4(K)+(\gamma_2(K)-2)(\epsilon-\lam_2)+\lam_3-\lam_2+\varsigma
\end{align*}
By Theorem \ref{thm:harper}, for any $K$ with prime decomposition number $p$ we have $\gamma_2(K) = 1,$ and $\gamma_3(K) = p+1$. Note that $\gamma_3(K) = p+1$ implies $\lambda_1 = p$. Direct substitution of these equalities yields
\begin{align*}
    &\gamma_4(K) = p+\epsilon-\lam_2\\
    &\gamma_5(K) = \binom{p}{2} + p+\lam_3-\lam_2+\varsigma
\end{align*}
which completes the proof.
\end{proof}
\section{A characterization of the \texorpdfstring{$\mathfrak{sl}_3$}{sl3} invariant on \texorpdfstring{$3$}{3}-braids}{\label{sec:braids}}
In this section we prove a formula which expresses the $\mathfrak{sl}_3$ invariant of an arbitrary 3-braid as a sum indexed by the size of twist regions and whose terms are functions of the set of nonzero indices. We develop this formulation by expressing each braid as a power of the Garside element times a positive braid. We give a closed formula for powers of the Garside element as a linear combination of non-elliptic basis elements in the associated web algebra. The $\mathfrak{sl}_3$ invariant can then be expressed as a linear combination of the invariants associated to the closures of the positive braid multiplied with each basis vector. 

\subsection{The Garside element in the web algebra}
We begin with our conventions for the braid group and introduce the web algebra. We then give a formula for powers of the Garside element in terms of a natural basis of the web algebra.\\

Let $\sigma_1,\dots, \sigma_{n-1}$ denote the standard generators of the braid group $B_n$ where each $\sigma_i$ is defined by the positive twist of strands $i$ and $i+1$. A braid is \emph{positive} if every crossing in the braid is positive. Let $B_n^+\subset B_n$ denote the \emph{positive braid monoid} of positive braid words. The \emph{Garside element} of the braid group $\Delta_n\in B_n^+$ is defined by 
\[
\Delta_n=\sigma_1\sigma_2\sigma_1\cdots \sigma_{n-1}\cdots\sigma_1\,.
\]
For $n=3$, we illustrate the Garside element and its square.

\begin{figure}[ht!]
    \centering
    \begin{tikzpicture}
    \draw (0,0) -- (0,0.25);
    \draw (0.5,0) -- (0.5,0.25);
    \draw (1,0) -- (1,0.25);
    \begin{scope}[shift={(0,0.25)}]
        \draw (0.5,0) -- (0,0.5);
        \draw[over] (0,0) -- (0.5,0.5);
        \draw (1,0) -- (1,0.5);
    \end{scope}
    \begin{scope}[shift={(0,0.75)}]
        \draw (1,0) -- (0.5,0.5);
        \draw[over] (0.5,0) -- (1,0.5);
        \draw (0,0) -- (0,0.5);
    \end{scope}
    \begin{scope}[shift={(0,1.25)}]
        \draw (0.5,0) -- (0,0.5);
        \draw[over] (0,0) -- (0.5,0.5);
        \draw (1,0) -- (1,0.5);
    \end{scope}
    \begin{scope}[shift={(0,1.75)}]
        \draw (0,0) -- (0,0.25);
        \draw (0.5,0) -- (0.5,0.25);
        \draw (1,0) -- (1,0.25);   
    \end{scope}
    \end{tikzpicture}
    \quad\quad\quad\quad\quad\quad\quad\quad\quad\quad
    \begin{tikzpicture}
    \draw (0,0) -- (0,0.25);
    \draw (0.5,0) -- (0.5,0.25);
    \draw (1,0) -- (1,0.25);
    \begin{scope}[shift={(0,0.25)}]
        \draw (0.5,0) -- (0,0.5);
        \draw[over] (0,0) -- (0.5,0.5);
        \draw (1,0) -- (1,0.5);
    \end{scope}
    \begin{scope}[shift={(0,0.75)}]
        \draw (1,0) -- (0.5,0.5);
        \draw[over] (0.5,0) -- (1,0.5);
        \draw (0,0) -- (0,0.5);
    \end{scope}
    \begin{scope}[shift={(0,1.25)}]
        \draw (0.5,0) -- (0,0.5);
        \draw[over] (0,0) -- (0.5,0.5);
        \draw (1,0) -- (1,0.5);
    \end{scope}
    \begin{scope}[shift={(0,1.75)}]
        \draw (0,0) -- (0,0.25);
        \draw (0.5,0) -- (0.5,0.25);
        \draw (1,0) -- (1,0.25);   
    \end{scope}
    \begin{scope}[shift={(0,2)}]
        \draw (0.5,0) -- (0,0.5);
        \draw[over] (0,0) -- (0.5,0.5);
        \draw (1,0) -- (1,0.5);
    \end{scope}
    \begin{scope}[shift={(0,2.5)}]
        \draw (1,0) -- (0.5,0.5);
        \draw[over] (0.5,0) -- (1,0.5);
        \draw (0,0) -- (0,0.5);
    \end{scope}
    \begin{scope}[shift={(0,3)}]
        \draw (0.5,0) -- (0,0.5);
        \draw[over] (0,0) -- (0.5,0.5);
        \draw (1,0) -- (1,0.5);
    \end{scope}
    \begin{scope}[shift={(0,3.5)}]
        \draw (0,0) -- (0,0.25);
        \draw (0.5,0) -- (0.5,0.25);
        \draw (1,0) -- (1,0.25);   
    \end{scope}
    \end{tikzpicture}
    \caption{$\Delta_3$ and $\Delta_3^2$}
    \label{fig:braids}
\end{figure}
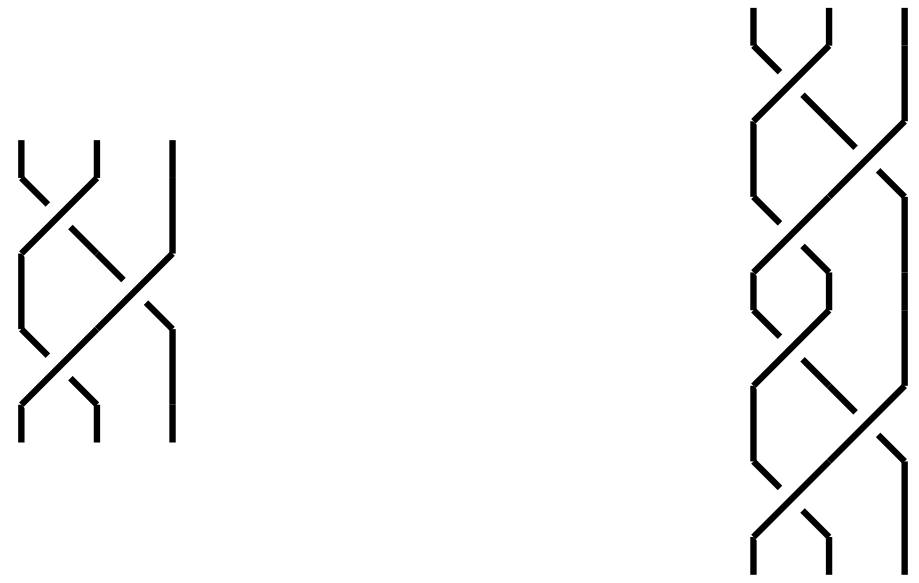
\begin{theorem}[{\cite[Theorems 5 and 7]{Gar}}]\label{thm:central}
    The center of $B_n$ is generated by $\Delta_n^2$. Every braid $b\in B_n$ can be written as a product $p\cdot \Delta_n^{m}$ for $p\in B_n^+$ and $m\in 2\mathbb Z$.
\end{theorem}

Let us also introduce the web algebra on three strands.

\begin{figure}[ht!]
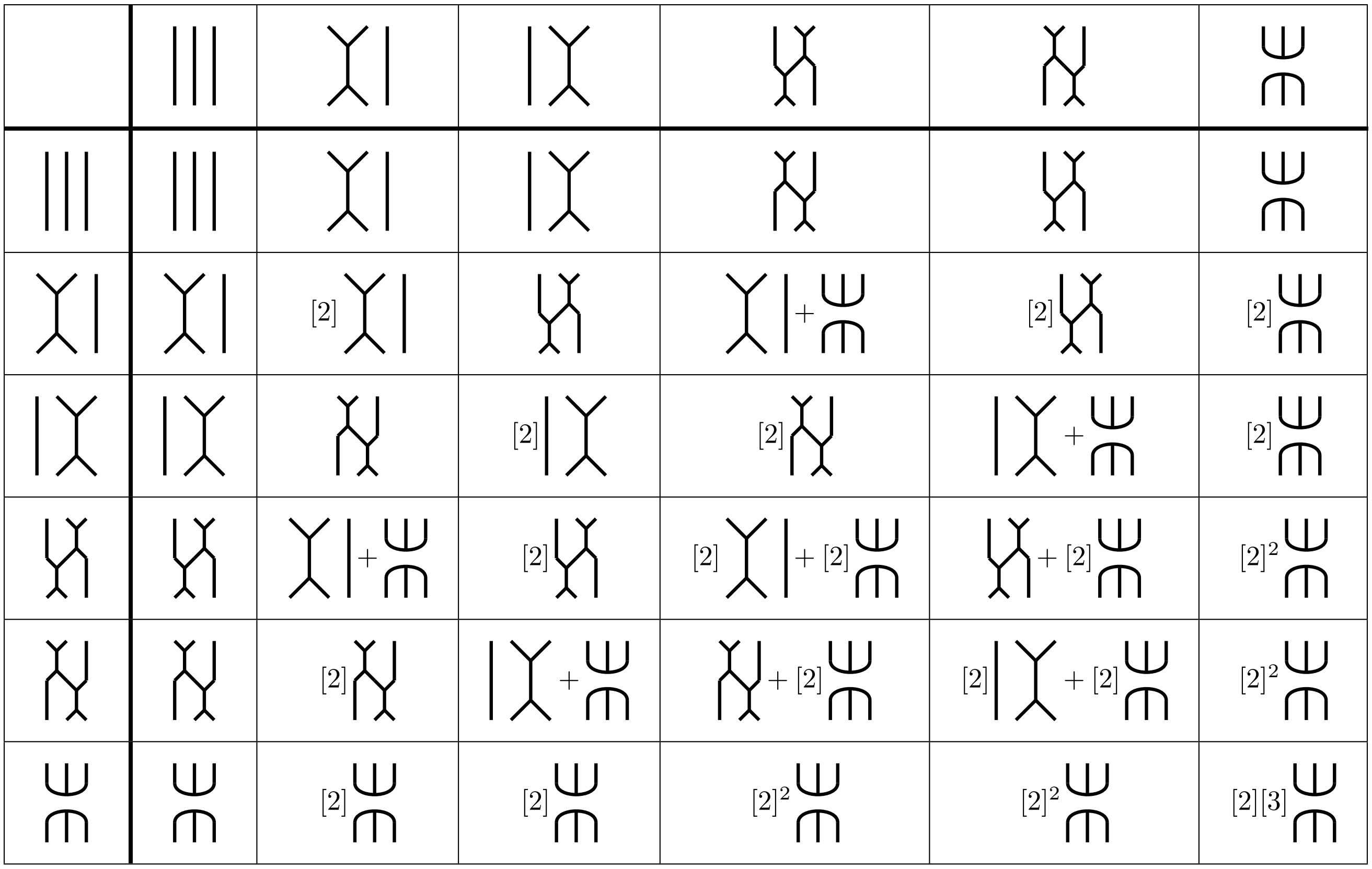
\label{multiplication_table}
    \centering
    \makebox[\textwidth]{$
    \renewcommand{\arraystretch}{1.2}
    \setlength{\arraycolsep}{11pt}
    \begin{array}{|c!{\vrule width 1.5pt}c|c|c|c|c|c|}
    \hline
    & {\ID} & {\LW} & {\RW} & {\LR} & {\RL} & {\CC} \\ \noalign{\hrule height 1.5pt}
    {\ID} & \ID & \LW & \RW & \RL & \LR & \CC \\ \hline
    {\LW} & \LW & [2]\,\LW & \LR & \LW + \CC & [2]\,\LR & [2]\,\CC \\ \hline
        {\RW} & \RW & \RL & [2]\,\RW & [2]\,\RL & \RW + \CC & [2]\,\CC \\ \hline

    {\LR} & \LR & \LW + \CC & [2]\,\LR & [2]\,\LW + [2]\,\CC & \LR + [2]\,\CC & [2]^2\,\CC \\ \hline
    {\RL} & \RL & [2]\,\RL & \RW + \CC & \RL + [2]\,\CC & [2]\,\RW + [2]\,\CC & [2]^2\,\CC \\ \hline
    {\CC} & \CC & [2]\,\CC & [2]\,\CC & [2]^2\,\CC & [2]^2\,\CC & [2][3]\,\CC \\ \hline
    \end{array}
    $}
    \caption{Multiplication table of non-elliptic basis elements of $\webalg$. The column label vector is stacked over the row label vector in the multiplication.}
    \label{fig:mult}
\end{figure}

\begin{definition}\label{web_arrow}
    Let $\mathrm{Web}_3$ be the $\mathbb Z[q, q^{-1}]$-algebra with a basis of ``non-elliptic" vectors
    \[\quad\idL\,, \quad\lw\,,\quad\rw\, ,\quad\lr\,,\quad\rl\,,  \quad\cc\,,\]
    which are labeled $v_0, \cdots, v_5$, multiplication in $\webalg$ given by stacking, and relations given by Equation \eqref{eq.webrelations}. An explicit multiplication table of basis elements is given in Figure \ref{fig:mult}.
\end{definition}

There is a natural quotient map $B_3\to \mathrm{Web}_3$, which applies the resolution \eqref{eq.crossingrelations} to each crossing. The diagram may then be simplified to the non-elliptic basis by applying the relations \eqref{eq.webrelations}. We denote this map from braids to webs by $\llangle\cdot\rrangle$, as it is the natural extension of the $\mathfrak{sl}_3$ invariant on link diagrams to  braids. In this way,
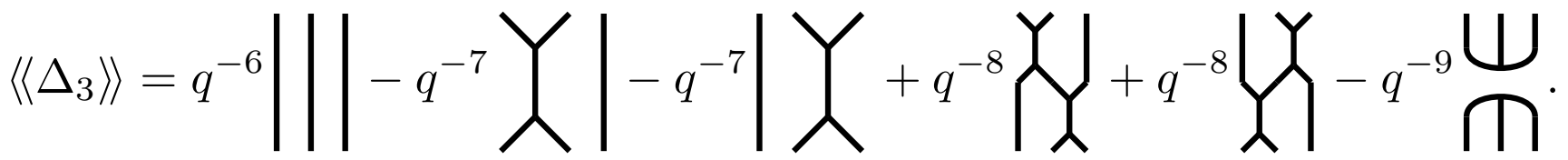
\begin{equation}\label{basis_garside}
    \llangle \Delta_3\rrangle =
    q^{-6}\,\idL\, 
    -q^{-7}\,\lw\,
    -q^{-7}\,\rw\, 
    +q^{-8}\,\rl\,
    +q^{-8}\,\lr\,
    -q^{-9}\,\cc\,. 
\end{equation}

We give a formula for the components of powers of the Garside element. First, we write 
\begin{align*}
    \llangle \Delta_3^n\rrangle&=g_0^{(n)}v_0+g_1^{(n)}v_1+g_2^{(n)}v_2+g_3^{(n)}v_3+g_4^{(n)}v_4+g_5^{(n)}v_5
    \\&=a_0^{(n)}v_0+a_1^{(n)}v_1+a_1^{(n)}v_2+a_2^{(n)}v_3+a_2^{(n)}v_4+a_3^{(n)}v_5\,.
\end{align*}
Notice that the coefficients on the pairs $v_1$ and $v_2$, and $v_3$ and $v_4$ are equal because $\Delta_3$ is invariant under reflection over a vertical axis.
Also $g_i^{(n)}$ and $a_i^{(n)}$ do not refer to taking $n^{th}$ powers, but just refer to the coefficient for each basis element $v_j$ of $\Delta_3^n$ for $0\leq j \leq 5$ which are not necessarily powers of $g_i^{(1)}$ and $a_i^{(1)}$.\\

In order to find the different powers of the polynomial linked to the Garside element in $B_3$, we just need to find the powers of different $a_i$ for $0\leq i \leq 3$. As previously explained, we know that we can resolve  the open Garside element $\Delta_3$ as a web. We begin by providing recursive formulas for these coefficients.

 \begin{lemma}\label{lem:recursive}
     The basis components of powers of the Garside element are given by the following recursive formulas:\begin{equation*}
    \begin{aligned}
    a_0^{(n+1)}&= q^{-6(n+1)}\\
    a_1^{(n+1)}&=-q^{-6n-7}+q^{-9}\cdot a_2^{(n)}\\
    a_2^{(n+1)}&=q^{-6n-8}+q^{-9}\cdot a_1^{(n)}\\
    a_3^{(n+1)}&= -q^{-6n-9}-2q^{-10}\cdot a_1^{(n)} -2q^{-11}\cdot a_2^{(n)}-q^{-12}\cdot a_3^{(n)}
    \end{aligned}
\end{equation*}
with $a_i^{(0)}=\delta_{i,0}$.
\end{lemma}
\begin{proof}
The base case of the induction is $\Delta_3^{(0)}=v_0$. 
    We obtain the recursion by computing the components of $\Delta_3^{n+1}$ in terms of the components of $\Delta_3^{n}$ and $ \Delta_3$,
\begin{align*}
\Delta_3^{n+1}&= \Delta_3^{n}\cdot \Delta_3\\ &=(a_0^{(n)}v_0+a_1^{(n)}v_1+a_1^{(n)}v_2+a_2^{(n)}v_3+a_2^{(n)}v_4+a_3^{(n)}v_5)
\\&\qquad\cdot (a_0^{(1)}v_0+a_1^{(1)}v_1+a_1^{(1)}v_2+a_2^{(1)}v_3+a_2^{(1)}v_4+a_3^{(1)}v_5)
\end{align*}

Figure \ref{fig:mult} tells us the combinations inside this multiplication that will contribute to each coefficient. Firstly, we observe that the only way to contribute to the coefficient of $a_0^{(n+1)}$ is by stacking $a_0^{(1)}$ on top of itself over and over. Hence, the recursive formula for the first coefficient is given by just multiplying the coefficient of the identity basis element, i.e. $a_0^{(1)}=q^{-6}$ in Equation \eqref{basis_garside},
\begin{align*}
    a_0^{(n+1)}&= q^{-6n}\cdot q^{-6}=q^{-6(n+1)}.
\end{align*}
The rest of the coefficients require more involved computations, but all follow a very similar procedure. We compute $a_1^{(n+1)}$ as an example:
\begin{align*}
    a_1^{(n+1)}&= a_1^{(n)}a_0^{(1)}+ a_0^{(n)}a_1^{(1)}+ [2]a_1^{(n)}a_1^{(1)}+ a_2^{(n)}a_1^{(1)}+a_1^{(n)}a_2^{(1)}+ [2]a_2^{(n)}a_2^{(1)}\\
    &= a_1^{(n)}(q^{-6})+ a_0^{(n)}(-q^{-7})+ [2]a_1^{(n)}(-q^{-7})+ a_2^{(n)}(-q^{-7})+a_1^{(n)}(q^{-8})+ [2]a_2^{(n)}(q^{-8})\\
    &= -q^{-6n-7} + q^{-9}\cdot a_2^{(n)}.\qedhere
\end{align*}
\end{proof}

Observe that every two steps in recursion for the coefficients $a_1^{(n)}$ (resp. $a_2^{(n)}$) is expressed in terms of $a_1^{(n-2)}$ (resp. $a_2^{(n-2)}$). Thus, we give closed formulas for the coefficients for the even powers of the Garside element. For our purposes, it suffices to consider $\Delta_3^m$ for $m\in 2\mathbb{Z}$ as in Proposition \ref{thm:central}.

\begin{theorem}
    The basis components of even powers $m=2k\geq 0$ of the Garside element are given by the following closed formulas:  
    \begin{align*}
    a_0^{(2k)}= q^{-12k},\qquad 
     a_1^{(2k)}=-q^{-15k}\{2\}[2][k]_{q^3}\,,\qquad
    a_2^{(2k)}= q^{-15k}\{2\}[k]_{q^3}\,,\qquad
    a_3^{(2k)}= -q^{-18k}\{2\}\{3\}[k]_{q^6}\,,
\end{align*}
where $\{2\}=q-q^{-1}$ and $\{3\}=q^2-1+q^{-2}$.
\end{theorem}
\begin{proof}
The formula for $a_0^{(m)}$ was proven in Lemma \ref{lem:recursive}.
    We will prove the remaining equalities by mathematical induction using the recursive formulas of Lemma \ref{lem:recursive}. The base case $m=0$ holds for these equalities as $m=k=0$ and the quantum integers $[k]_{q^3}$ and $[k]_{q^6}$ vanish.\\

Inductive Step: Let us assume:
\[
a_1^{(m-2)}=-q^{-15(k-1)}\{2\}[2][k-1]_{q^3}\,,
\quad
a_2^{(m-2)}=q^{-15(k-1)}\{2\}[k-1]_{q^3}\,,
\quad 
a_3^{(m-2)}=-q^{-18(k-1)}\{2\}\{3\}[k-1]_{q^6}\,.
\]
Now, we want to prove our formulae for $a_i^{(m)}$. By applying the recursion twice,
\begin{align*}
    a_1^{(m)}&=-q^{-6m-1}+q^{-9}\cdot a_2^{(m-1)}
    =
    -q^{-6m-1}+q^{-9} 
    (q^{-6(m-1)-2}+q^{-9}\cdot a_1^{(m-2)})
    \\&=
    -q^{-12k-1}+q^{-12k-5} - q^{-15k-3}\{2\}[2][k-1]_{q^3}
    \\&=
    -q^{-15k}\left(q^{3k-3}\{2\}[2]+q^{-3}\{2\}[2][k-1]_{q^3}
    \right)
    \\&=
    -q^{-15k}\{2\}[2]\left(q^{3(k-1)}+q^{-3}[k-1]_{q^3}
    \right)
    \\&= -q^{-15k}\{2\}[2][k]_{q^3}
\end{align*}
which completes the induction for $a_1^{(m)}$. The argument for the others is similar.
\end{proof}

\begin{corollary}\label{cor.invert}
    For $m=2k\leq 0$, the components of the Garside elements are given by 
    \begin{align*}
    a_0^{(2k)}= q^{-12k},\qquad 
     a_1^{(2k)}=q^{-15k}\{2\}[2][k]_{q^3}\,,\qquad
    a_2^{(2k)}= -q^{-15k}\{2\}[k]_{q^3}\,,\qquad
    a_3^{(2k)}= q^{-18k}\{2\}\{3\}[k]_{q^6}\,.
\end{align*}
\end{corollary}
\begin{proof}
    The inversion map on crossings induces the map $q\mapsto q^{-1}$ on webs by \eqref{eq.crossingrelations}. Thus, the components of negative powers of the Garside element are given by inverting $q$. In particular, $\{2\}\mapsto -\{2\}$ provides additional signs on $a_1^{(m)}$, $a_2^{(m)}$, and $a_3^{(m)}$.
\end{proof}
\subsection{\texorpdfstring{$\mathfrak{sl}_3$}{sl3} invariant on positive \texorpdfstring{$3$}{3}-braids}
{To use the result that every braid can be represented as a negative even power of the Garside element and a positive braid, we now characterize the invariant on positive braids. We utilize the linearity of the trace map to express the evaluation of a positive braid in terms of an algorithm, which inputs information from each state through a binary sequence.}
\begin{lemma}
    Let $p\in B_3^+$ considered up to Markov I moves. Then $p$ may be written in the form
    \[
    \sigma_1^{0}\sigma_2^{0}, \qquad \sigma_1^{n_1}\sigma_2^{0},\qquad 
    \sigma_1^0\sigma_2^{n_2},\qquad 
    \text{or} \qquad \sigma_1^{n_1}\sigma_2^{n_2}\cdots \sigma_2^{n_d}
    \]
    where each $n_i>0$.
\end{lemma}
\begin{proof}
    Since $p$ is a positive braid, there cannot be any cancellation among its generators. Thus, the split braids $p=e$, $p=\sigma_1^{n_1}$, and $p=\sigma_2^{n_2}$ are the only presentations of their underlying braid class up to Markov I moves. This gives the first two presentations of $p$.\\\\
    We now consider non-split braids. Without loss of generality, we can assume $d$ is even. In other words, the first and last letters of $p$ are different. If it were odd, we apply the Markov I move to simply slide the last letter $\sigma_k$ in the braid word to before the first letter of the word which is also $\sigma_k$.
    If the positive braid word starts with $\sigma_2$, we just slide it down to the end to start on $\sigma_1$ instead. We can assume every $n_i>0$, since if $n_i=0$, we could just merge the two $\sigma$s of the same index before and after it.
\end{proof}

We define the \emph{trace map} $\mathrm{tr}: \webalg \to \mathbb Z[q,q^{-1}]$, which is $\mathbb Z[q, q^{-1}]$-linear, that represents the evaluation of an open web as a closed one. \\\\
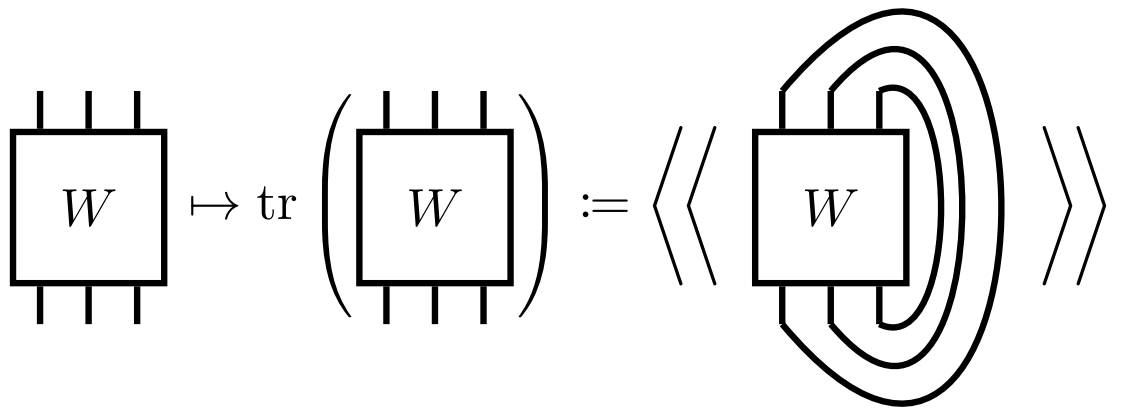
\begin{figure}[H]
    \centering
    \vspace{-80pt}
    \begin{tikzpicture}
        \node (W) [draw, minimum size=1cm] {$W$};
        \draw (W.north west) ++(0.2,0) -- ++(0,0.25);
        \draw (W.north)            -- ++(0,0.25);
        \draw (W.north east) ++(-0.2,0) -- ++(0,+0.25);
        \draw (W.south west) ++(0.2,0) -- ++(0,-0.25);
        \draw (W.south)            -- ++(0,-0.25);
        \draw (W.south east) ++(-0.2,0) -- ++(0,-0.25);
    \end{tikzpicture}
    $\mapsto \mathrm{tr}\left(
    \begin{tikzpicture}
        \node (W) [draw, minimum size=1cm] {$W$};
        \draw (W.north west) ++(0.2,0) -- ++(0,0.25);
        \draw (W.north)            -- ++(0,0.25);
        \draw (W.north east) ++(-0.2,0) -- ++(0,+0.25);
        \draw (W.south west) ++(0.2,0) -- ++(0,-0.25);
        \draw (W.south)            -- ++(0,-0.25);
        \draw (W.south east) ++(-0.2,0) -- ++(0,-0.25);
    \end{tikzpicture}
    \right)\coloneqq\leftlangle\,\,\,
    \begin{tikzpicture}
        \node (W) [draw, minimum size=1cm] {$W$};
        \coordinate (T1) at ($(W.north west)+( 0.2, 0.25)$);
        \coordinate (T2) at ($(W.north)     +( 0,   0.25)$);
        \coordinate (T3) at ($(W.north east)+(-0.2, 0.25)$);
        \coordinate (B1) at ($(W.south west)+( 0.2,-0.25)$);
        \coordinate (B2) at ($(W.south)     +( 0,  -0.25)$);
        \coordinate (B3) at ($(W.south east)+(-0.2,-0.25)$);
        \draw (W.north west) ++(0.2,0)  -- (T1);
        \draw (W.north)                 -- (T2);
        \draw (W.north east) ++(-0.2,0) -- (T3);
        \draw (W.south west) ++(0.2,0)  -- (B1);
        \draw (W.south)                 -- (B2);
        \draw (W.south east) ++(-0.2,0) -- (B3);
        \draw (B3) to[out=-25, in=25] (T3);
        \draw (B2) to[out=-50, in=50, looseness = 3] (T2);
        \draw (B1) to[out=-50, in=50, looseness = 5] (T1);
    \end{tikzpicture}
    \hspace{-16pt}
    \,\,\,\,\,\rightrangle$
    \vspace{-40pt}
    \caption{The trace of an element $W$ in $\webalg$}
    \label{fig:trace}
\end{figure}
This map also fits the commutative diagram
\[
\begin{tikzcd}
    B_3 \arrow[r,"\llangle\cdot\rrangle"] \arrow[d,"\widehat{\cdot}"'] & \webalg \arrow[d,"\mathrm{tr}"]
    \\
    \mathrm{Links} \arrow[r,"\llangle\cdot\rrangle"]& \mathbb{Z}[q,q^{-1}]
\end{tikzcd}\,.
\]

Thus, for a braid $b$ where $\llangle b \rrangle=\sum_{i=0}^5 b_i v_i$ for some $b_i\in \mathbb{Z}[q,q^{-1}]$, linearity of the trace implies
\begin{equation}\label{eq.tracecomponents}
    \llangle \,\widehat{b}\,\rrangle =\sum_{i=0}^5 b_i \llangle \widehat{v_i}\rrangle.
\end{equation}

Let $0\le k_i\le n_i$ for all $1\le i\le d$. Then, denote $\sum_{i=1}^dk_i=\mathbf{k}$, and $\sum_{i=1}^dn_i=\mathbf{n}$. \\\\
In a slight revision to our earlier notation, we use $k_i$ to indicate the number of $W$-resolutions in the $i$-th twist region. Notice that every resolution of a positive 3-braid $p$ is uniquely determined by the $n$-tuple $(k_1, \dots, k_d)$ indicating the number of $W$-resolutions in each twist region. Let $\mathbb W((k_1, \dots, k_d))$ be the associated web diagram and $\widehat{\mathbb W}((k_1, \dots, k_d))$ its closure.
Now, we want to compute the $\mathfrak{sl}_3$ invariant of the braid closure of $p$:
\begin{align*}
    \llangle \,\widehat{p}\,\rrangle&=\sum_{k_1, \dots, k_d}\binom{n_1}{k_1}\cdots \binom{n_d}{k_d}\phi((k_1, \dots, k_d))\llangle \widehat{\mathbb W}((k_1, \dots, k_d))\rrangle 
    \\&=
\sum_{k_1, \dots, k_d}\binom{n_1}{k_1}\cdots \binom{n_d}{k_d}(-q)^{-2\mathbf{n}-\mathbf{k}} \llangle \widehat{\mathbb W}(k_1, \dots, k_d)\rrangle
\end{align*}
after computing $\phi((k_1, \dots, k_d))=(-q)^{-2\mathbf{n}-\mathbf{k}}$. Next, we pull out any bubbles for $k_i\ge 2$. For each $i$ and given $k_i$, define \[t_i = \begin{cases}
    1 & : k_i \ge 1 \\
    0 & : \text{else}
\end{cases}.\]
This now reduces our consideration to webs indexed by cyclic binary sequences.

\begin{lemma}\label{lem:sl3-positive-3-braids}
Given $p=\sigma_1^{n_1}\sigma_2^{n_2}\cdots \sigma_2^{n_d}\in B_3^+$ with each $n_i>0$,
\[
\llangle \,\widehat{p}\,\rrangle = \sum_{k_1, \dots, k_d}\binom{n_1}{k_1}\cdots \binom{n_d}{k_d}(-q)^{-2\mathbf{n}-\mathbf{k}} [2]^{\sum_{k_i: k_i\ge 2}(k_i-1)}\widehat{\mathbb W}((t_1, \dots, t_d))\rrangle.
\]
\end{lemma}

We describe an algorithm to compute $\llangle\widehat{\mathbb W}((t_1, \dots, t_d))\rrangle$ on even length binary sequences $(t_1, \dots, t_d)$. Note that the evaluations of (crossingless) webs are invariant under taking mirror image. This implies that isotoping the $d$-th twist region produces a web $\mathbb W((t_d, t_1 \dots, t_{d-1}))\rrangle$ which has the same evaluation as $\mathbb W((t_1, \dots, t_d))\rrangle$. Hence $(t_1, \dots, t_d)$ may be considered as a cyclic binary sequence.
\\\\
Define a function $w$ on non-negative even integers
\[
w(2r)=
\begin{cases}
[3]\big([2]^2 + \sum_{j=1}^{r-1} [2]^{2j}\big) &:r>0\\
[3]^2 &:r=0
\end{cases}\,.
\]

\begin{lemma}\label{lem:Wbin}
    Let $T=(t_1, \dots, t_d)$ be a cyclic binary sequence. Let $z$ be the number of maximal zero-blocks in $T$ with odd length and let $\mathbf{t}=\sum_{i=1}^d t_i$. Then
    \[
    \llangle \mathbb W(T)\rrangle=
    \begin{cases}
        [3]^3 &:T=(0,\dots, 0)
        \\
        [2]^zw(\mathbf{t}-z) &:T\neq(0,\dots, 0)
    \end{cases}\,.
    \]
\end{lemma}
\begin{proof}
    Let us describe an algorithm to evaluate $\llangle \mathbb W((t_1, \dots, t_d))\rrangle$ which will produce the desired formula.
\begin{enumerate}
    \item Delete any pairs of adjacent $00$s until no pairs are left. This is because 00 corresponds to two consecutive $O$ resolutions, which cancels out the parity, and is equivalent to the identity.
    
    \item If the remaining sequence is empty, $T(t_1, \cdots, t_d)=[3]^3$. This is because the resulting empty sequence implies $t_i=0$ for all $i$ in the original input sequence, which corresponds to an All-$O$ web. Closing the All-$O$ web gives us three unlinked unknots, which evaluates to $[3]^3$.

    \item If $T$ did not correspond to the All-$O$ resolution, the resulting even length cyclic sequence is now made up of runs of $1$s separated by isolated $0$s. We let $z$ be the number of $0$s, which is equivalent to the number of maximal zero-blocks in $T$.
    
    If $d=2$ and $z=1$, then the sequence is $01$, up to cycling. This evaluates to $[2][3]^2=[2]^zw(0)$.
    
    Otherwise, notice that each of the $z$ $0$s occurs in a sequence of the form $101$. Since the two $1$s beside it share the same coordinate parity, this gives us a $[2]$ bubble and reduces $101\to 1$. Pulling out $[2]^z$ from all $101$ sequences, reduces us to an even length sequence of all $1$s. More precisely, we end up with $\mathbf{t}-z$ total $1$s, which brings us back to the even squared situation covered in {\cite[Lemma 5.10]{HK}}, giving us $w(\mathbf{t}-z)$. 
\end{enumerate}
\end{proof}

\begin{lemma}\label{lem.Tstack}
    For each 3-web $\mathbb{W}(T)$ associated to a cyclic binary sequence $T=(t_1,\dots, t_d)$, the evaluation of the closure of $v_i\cdot \mathbb{W}(T)$ for $0\leq i\leq 4$ may be expressed in terms of the evaluated closure $\mathbb{W}(T'(i))$ where $T'(i)$ is the cyclic binary sequence     
\[
T'(i)=
\begin{cases}
    (t_1,\dots, t_d) & i=0
    \\
    (1,0,t_1,\dots, t_d) & i=1
    \\
    (0,1,t_1,\dots, t_d) & i=2
    \\
    (1,1,t_1,\dots, t_d) & i=3
    \\
    (t_1,\dots, t_d,0,1,1,0) & i=4.
\end{cases}
\]
\end{lemma}
\begin{proof}
    {We can write each of the first five basis vectors as short binary sequence of resolutions represented by $0$s and $1$s.
    \begin{itemize}
        \item Multiplying $v_0$, which is the identity web, to the bottom of a positive braid does not change it.
        \item Multiplying $v_1$ to the bottom of a positive braid would require a $1$ in an even index to indicate the $W$-resolution and then a $0$ to maintain the parity of the coordinates, giving us $(1,0,t_1, \cdots, t_d)$.
        \item For $v_2$, to maintain that the $1$ indicates a $W$-resolution on the right side, we add a $0$ to cushion the parity, then a $1$ to indicate. 
        \item The cases for $v_3$ and $v_4$ proceed in a similar manner. 
    \end{itemize}}
\end{proof}

Combining the the previous results gives an algorithm to evaluate any positive 3-braid. 

\begin{theorem}\label{thm:sl3-for-3-braids}
Take $b\in B_3$, with $b=p\cdot \Delta_3^{-2k}$ for some $k\in \mathbb Z_{\ge 0}$ and $p=\sigma_1^{n_1}\cdots \sigma_2^{n_d}$ up to Markov moves. Then, 
\begin{equation*}
    \begin{aligned}
        \llangle \hat{b}\rrangle &= q^{12k-2\mathbf{n}}\sum_{k_1, \cdots, k_d}\binom{n_1}{k_1}\cdots \binom{n_d}{k_d}((-q)^{-\mathbf{k}})[2]^{\sum_{k_i:k_i\ge 2}(k_i-1)}
        \\&\qquad
        \cdot 
        \bigg(\llangle \mathbb{W}(T'(0))\rrangle +q^{3k}\{2\}[2][k]_{q^3}\mathbb{W}(T'(1))\rrangle+q^{3k}\{2\}[2][k]_{q^3}\llangle \mathbb{W}(T'(2))\rrangle
        \\&\qquad\qquad -q^{3k}\{2\}[k]_{q^3}\mathbb{W}(T'(3))\rrangle
        -q^{3k}\{2\}[k]_{q^3}\mathbb{W}(T'(4))\rrangle+q^{6k}\{2\}\{3\}[k]_{q^3} [2]^{\mathbf{t}+1}[3]\bigg) 
    \end{aligned}
\end{equation*}
\end{theorem}
\begin{proof}
    Recall from \eqref{eq.tracecomponents} that linearity of the trace implies that
    \[
\llangle \,\widehat{b}\,\rrangle = \llangle \,\widehat{p\cdot \Delta_3^{-m}}\,\rrangle = \sum_{i=0}^5 g_i\llangle \widehat{p\cdot v_i}\rrangle.
\]
Our next step is to characterize $\llangle \widehat{p\cdot v_i}\rrangle$ for each of $0\le i\le 5$ in a way similar to $\llangle \,\widehat{p}\,\rrangle$. Each multiplication of a basis element $v_i$ with $p$ modifies the binary sequence function $\llangle \mathbb{W}(T)\rrangle$ appearing in Lemmas \ref{lem:sl3-positive-3-braids} and \ref{lem:Wbin} with  
\[
T = (t_1,\dots, t_d) \qquad \mbox{where}\qquad
t_i = \begin{cases}
    1 & : k_i \ge 1 \\
    0 & : \text{else}
\end{cases}.
\]
In each case 
\[
\llangle \widehat{p\cdot v_i}\rrangle= \sum_{k_1, \cdots, k_d}\binom{n_1}{k_1}\cdots \binom{n_d}{k_d}(-q)^{-2\mathbf{n}-\mathbf{k}} [2]^{\sum_{k_i: k_i\ge 2}(k_i-1)}\llangle \mathbb{W}(T'(i))\rrangle
\]
where $\mathbb{W}(T'(i))$ is defined as in Lemma \ref{lem.Tstack}. For the case $i=5$, the web evaluates exactly to $[2]^{\mathbf{t}+1}[3]$. Each elementary web multiplied with $v_5$ contributes a factor of $[2]$ and the closure of $v_5$ is $[2][3]$. Finally, summing over these terms gives the desired formula.
\end{proof}

This result has been verified by the code accessible at \cite{AngieCode}. The code allows for computation of the $\inv$ invariant of 3-braids in terms of the Garside exponent and positive braid word.

\bibliographystyle{alpha}
\bibliography{references}
\end{document}